\documentclass[12pt,leqno]{amsart}
\usepackage{amssymb}
\usepackage{amsmath,amssymb,color}
\newtheorem{thm}{Theorem}[section]
\newtheorem{lem}[thm]{Lemma}
\newtheorem{cor}[thm]{Corollary}
\newtheorem{prop}[thm]{Proposition}

\numberwithin{equation}{section}

\newcommand{\weight}{e^{2s\alpha}}

\newcommand{\ep}{\varepsilon}
\newcommand{\la}{\lambda}
\newcommand{\va}{\varphi}
\newcommand{\ppp}{\partial}

\newcommand{\www}{\widetilde}

\newcommand{\R}{\mathbb{R}}

\newcommand{\N}{\mathbb{N}}

\newcommand{\ooo}{\overline}
\newcommand{\OOO}{\Omega}

\newcommand{\sumij}{\sum_{i,j=1}^d}

\newcommand{\HHHH}{H^{2,1}}

\newcommand{\hhalf}{\frac{1}{2}}

\newcommand{\ddd}{\mathrm{div}\,}  

\allowdisplaybreaks

\title
[]
{
Lipschitz stability for determination of states and 
inverse source problem for the mean field game equations
}

\author{$^1$ Oleg Imanuvilov, $^2$ Hongyu Liu and $^3$ Masahiro Yamamoto}

\thanks{
$^1$ Department of Mathematics, Colorado State University, 101 Weber Building, 
Fort Collins CO 80523-1874, USA, e-mail: {\tt oleg@math.colostate.edu}
\\
$^2$ Department of Mathematics, City University of Hong Kong, Kowloon, 
Hong Kong SAR, China, email: {\tt hongyliu@cityu.edu.hk}
\\
$^3$ Graduate School of Mathematical Sciences, The University
of Tokyo, Komaba, Meguro, Tokyo 153-8914, Japan, 
e-mail: {\tt myama@ms.u-tokyo.ac.jp}
}

\date{}
\begin{document}
\maketitle

\begin{abstract}
In a bounded domain $\OOO \subset \R^d$ over time interval 
$(0,T)$, we consider mean field game equations 
whose principal coefficients depend on the time and state variables
with a general Hamiltonian.
We attach  the non-zero Robin boundary condition.
We first prove the Lipschitz stability in $\OOO \times (\ep, T-\ep)$
with given $\ep>0$ 
for the determination of the solutions by Dirichlet data on arbitrarily 
chosen subboundary of $\ppp\OOO$.
Next we prove the Lipschitz stability for an inverse problem of determining 
spatially varying factors of source terms and a coefficient
by extra boundary data and spatial data at an intermediate time.    
\end{abstract} 
\baselineskip 18pt

\section{Introduction}

Recently the mean field game has called great attention and we refer
for example, to
Achdou, Cardaliaguet, Delarue, Porretta and Santambrogio
\cite{ACDPS}, Cardaliaguet, Cirat and Porretta \cite{CCP},
Lasry and Lions \cite{LL} and the references therein.

Let $\OOO \subset \R^d$ be a smooth bounded domain and let 
$\nu = \nu(x) = (\nu_1(x), ..., \nu_d(x))$ be the outward unit normal 
vector to $\ppp\OOO$ at 
$x \in \ppp\OOO$.  Let  
$$
Q:= \OOO \times (0,T).
$$
Then a forward problem for one typical mean field game system 
can be described by 
\begin{equation}\label{(1.1)}
\left\{ \begin{array}{rl}
& \ppp_tu(x,t) + a(x,t)\Delta u(x,t) 
- \hhalf\kappa(x,t)\vert \nabla u(x,t)\vert^2 + c_0(x,t)v = F(x,t), \\
& \ppp_tv(x,t) - \Delta (a(x,t)v(x,t)) - \ddd(\kappa(x,t)v(x,t)\nabla u(x,t))
= G(x,t), \quad (x,t) \in Q
\end{array}\right.
                                                          \end{equation}
with the homogeneous Robin boundary condition
\begin{equation}\label{(1.2)}
\nabla u \cdot \nu = \nabla (av)\cdot \nu = 0 
\quad \mbox{on $\ppp\OOO \times (0,T)$},                          
                   \end{equation}
and given $u(\cdot,T)$ and $v(\cdot,0)$ in $\OOO$.

Although there are many works on the well-posedness and related properties 
of solutions $u, v$ for the above system,
very few efforts have been devoted to inverse problems for (\ref{(1.1)}).
We can refer to  Klibanov \cite{Kl23}, Klibanov and Averboukh \cite{KlAv},
Klibanov, Li and Liu \cite{KLL1}, \cite{KLL2},
Liu and Zhang \cite{LZ1}, \cite{LZ2}.

For inverse problems, we can mention two types for example:
\\
{\bf (i) Determination of state:} {\it Determine $u, v$ by extra data on a 
subboundary 
of $\ppp\OOO$ over a time interval.}
\\
{\bf (ii) Inverse source problem:} {\it Determine spatially varying 
factors of the source terms $F(x,t)$ and $G(x,t)$ by lateral data and spatial
data of $u, v$.}
\\

As for the type (i), we can refer to \cite{KLL2}, and as for other types of
state determination for $u, v$ with data chosen among
$u(\cdot,0), u(\cdot,T), v(\cdot,0), v(\cdot,T)$, see 
\cite{Kl23}, \cite{KlAv}, \cite{KLL1}.  Also see Liu and Yamamoto \cite{LY} 
as related state determination problems.  We can further mention the 
unique continuation as other important property 
for partial differential equations, and can refer to 
Liu, Imanuvilov and Yamamoto \cite{ILY} for unique continuation of solutions to  mean field game equations.

However, to the best knowledge of the authors, there are no publications on
the uniqueness and the stability for the inverse source  problem for the mean field game equations.
The main purpose of this article is to establish the Lipschitz stability
for the above two types of inverse problems, which have not been found in 
the existing articles.

We mainly consider a linearized equation of (\ref{(1.1)}), which is formulated 
as follows.

We set 
$\ppp_i = \frac{\ppp}{\ppp x_i}$, $1\le i \le d$ and 
$\ppp_t = \frac{\ppp}{\ppp t}$ and for $\gamma := (\gamma_1, ...,
\gamma_d) \in (\N \cup \{ 0\})^d$ we define 
$\ppp_x^{\gamma}: = \ppp_1^{\gamma_1}\cdots \ppp_d^{\gamma_d}$ and
$\vert \gamma\vert := \gamma_1 + \cdots + \gamma_d$,
and introduce the following functional spaces
\begin{align*}
& C^{1,0}(\ooo{Q}):= \{ u\in C(\ooo{Q});\, 
\nabla u \in C(\ooo{Q})\}, \quad
  C^1(\ooo{Q}):= \{ u\in C(\ooo{Q});\,
\nabla u, \ppp_tu \in C(\ooo{Q})\},\\
& C^{2,1}(\ooo{Q}):= \{ u\in C(\ooo{Q});\, 
\nabla u, \ppp_i\ppp_ju, \ppp_tu  \in C(\ooo{Q}), \, \, 
1\le i,j \le d \}, \\
& H^{2,1}(Q):= \{ u\in L^2(Q);\, \nabla u, \ppp_i\ppp_ju, \ppp_tu  
\in L^2(Q), \, \, 1\le i,j \le d \} \quad \mbox{with the norm}: \\
& \Vert u\Vert_{H^{2,1}(Q)}:= \left( \sum_{\vert \gamma\vert \le 2}
\Vert \ppp_x^{\gamma}u\Vert^2_{L^2(Q)}
+ \Vert \ppp_tu\Vert^2_{L^2(Q)}\right)^{\frac{1}{2}}.
\end{align*}

Throughout this article, we consider the following partial differential 
operators of the second order with $(x,t)$- dependent coefficients defined 
by
$$
\left\{ \begin{array}{rl}
&A(t)u:= \sumij a_{ij}(x,t)\ppp_i\ppp_j u + \sum_{i=1}^d a_j(x,t)\ppp_ju
+ a_0(x,t)u, \\
&B(t)v:= \sumij b_{ij}(x,t)\ppp_i\ppp_j v + \sum_{i=1}^d b_j(x,t)\ppp_jv
+ b_0(x,t)v,\\
& A_0(t)u := \sum_{\vert \gamma\vert \le 2}\widetilde  a_{\gamma}(x,t) 
\ppp_x^{\gamma}u,
\end{array}\right.
$$
where 
\begin{equation}\label{lopukh1}
\left\{ \begin{array}{rl}
&a_{ij}, b_{ij} \in C^1(\ooo{Q}), \quad a_{ij} = a_{ji},\,\,
b_{ij} = b_{ji}\quad \mbox{ for}\quad 1\le i,j\le d, \cr\\
& a_k, b_k \in L^{\infty}(Q)\quad \mbox{ for}\quad  0\le k \le d,\quad 
\widetilde a_{\gamma}\in L^{\infty}(Q)\quad  \mbox{for} \quad 
\vert\gamma \vert \le 2,
\end{array}\right.
\end{equation}
and it is assumed that there exists a constant 
$\chi > 0$ such that
\begin{equation}\label{lopukh}
\sumij a_{ij}(x,t)\xi_i\xi_j \ge \chi \sum_{j=1}^d \xi_j^2 \quad\mbox{and}\quad 
\sumij b_{ij}(x,t)\xi_i\xi_j \ge \chi \sum_{j=1}^d \xi_j^2
\end{equation}
for all $(x,t) \in Q$ and $\xi_1, ..., \xi_d \in \R$.

We define
$$
\ppp_{\nu_A}u := \sumij a_{ij}(\ppp_ju)\nu_i, \quad
\ppp_{\nu_B}v := \sumij b_{ij}(\ppp_jv)\nu_i  \quad
\mbox{on $\ppp\OOO\times (0,T)$}.
$$

We consider a linearized mean field game system:
\begin{equation}\label{(1.3)}
\left\{ \begin{array}{rl}
& \ppp_tu + A(t)u = c_0(x,t)v + F(x,t), \\
& \ppp_tv - B(t)v = A_0(t)u + G(x,t), \quad (x,t) \in Q
\end{array}\right.
\end{equation}
with the Robin boundary condition
\begin{equation}\label{(1.4)}
\left\{ \begin{array}{rl}
& \ppp_{\nu_A}u(x,t) - p(x,t)u(x,t) = g(x,t), \\ 
& \ppp_{\nu_B}v(x,t) - q(x,t)v(x,t) = h(x,t), \quad (x,t) \in \ppp\OOO \times 
(0,T),
\end{array}\right. 
\end{equation}
where we assume
\begin{equation}\label{1.7}
p, q \in C^1(\ppp\OOO\times [0,T]) \quad \mbox{and}\quad c_0\in L^\infty(Q).
\end{equation}

We emphasize that general Hamiltonians also can be considered. Such Hamiltonians produce  
the second-order partial differential operator 
$A_0(t)u$ after the linearization.

In this article, we establish the Lipschitz stability results
for the following inverse problems.  Let $\Gamma$ be an arbitrarily chosen 
non-empty subboundary of $\ppp\OOO$, $t_0 \in (0,T)$ be arbitrarily given, and
functions $u, v$ satisfy (\ref{(1.3)}) and (\ref{(1.4)}). We consider the following inverse problems
\\
{\bf Determination of state:}
{\it Determine $u, v$ in $Q$ by $u\vert_{\Gamma\times (0,T)}$ and
$v\vert_{\Gamma\times (0,T)}$.
}
\\
{\bf Inverse source problem:}
{\it
In (\ref{(1.3)}), let $F(x,t) = q_1(x,t)f_1(x)$ and $G(x,t) = q_2(x,t)f_2(x)$ 
for $(x,t) \in Q$ where functions $q_1$ and $q_2$ are given.  
Determine $f_1, f_2$ in $\OOO$ by data
$u\vert_{\Gamma\times (0,T)}$, $v\vert_{\Gamma\times (0,T)}$
and $\{ u(\cdot,t_0), \, v(\cdot,t_0)\}$ in $\OOO$.
}
\\

We first state our main result for the state determination.
\\
\begin{thm}\label{t1}
{\it 
We assume (\ref{lopukh1}), (\ref{lopukh}) and (\ref{1.7}).
Moreover let $u,v \in H^{2,1}(Q)$ satisfy (\ref{(1.3)}) and (\ref{(1.4)}).
For arbitrarily given $\ep > 0$, we can find a constant $C_{\ep} >0$ such that 
\begin{align*}
& \Vert u\Vert_{\HHHH(\OOO\times (\ep,T-\ep))}
+ \Vert v\Vert_{\HHHH(\OOO\times (\ep,T-\ep))}
\le C_{\ep}(\Vert F\Vert_{L^2(Q)} + \Vert G\Vert_{L^2(Q)}\\
+ & \Vert u\Vert_{H^1(\Gamma \times (0,T))}
+ \Vert v \Vert_{H^1(\Gamma \times (0,T))}
+ \Vert \ppp_tg\Vert_{L^2((\ppp\OOO\setminus \Gamma)\times (0,T))}
+ \Vert \ppp_th\Vert_{L^2((\ppp\OOO\setminus \Gamma)\times (0,T))}\\
+ & \Vert g\Vert_{L^2(0,T;H^{\hhalf}(\ppp\OOO))}
   + \Vert h\Vert_{L^2(0,T;H^{\hhalf}(\ppp\OOO))}).
\end{align*}
}
\end{thm}

In particular, we directly see
\begin{align*}
& \Vert u(\cdot,t)\Vert_{L^2(\OOO)}
+ \Vert v(\cdot,t)\Vert_{L^2(\OOO)}
\le C_{\ep}(\Vert F\Vert_{L^2(Q)} + \Vert G\Vert_{L^2(Q)}\\
+ & \Vert u\Vert_{H^1(\Gamma \times (0,T))}
+ \Vert v \Vert_{H^1(\Gamma \times (0,T))}
+ \Vert \ppp_tg\Vert_{L^2((\ppp\OOO\setminus \Gamma)\times (0,T))}
+ \Vert \ppp_th\Vert_{L^2((\ppp\OOO\setminus \Gamma)\times (0,T))}\\
+& \Vert g\Vert_{L^2(0,T;H^{\hhalf}(\ppp\OOO))}
 + \Vert h\Vert_{L^2(0,T;H^{\hhalf}(\ppp\OOO))})
\end{align*}
for $\ep \le t \le T-\ep$.

We emphasize that Theorem \ref{t1} asserts unconditional stability 
in the case of non-homogeneous Robin boundary condition. The 
unconditional stability means 
that we do not need to impose any boundedness assumptions for $u$ and $v$.
 
In Klibanov, Li and Liu \cite{KLL2}, 
the H\"older stability is proved with data $u, v, \nabla u, \nabla v$ 
on the whole lateral boundary $\ppp\OOO \times 
(0,T)$, which is the case of $\Gamma = \ppp\OOO$ in Theorem \ref{t1}.
By the parabolicity of the equations (\ref{(1.3)}), extra boundary 
data should be limited to any small subboundary.
\\

Next we state our main result on the inverse source problem.
In (\ref{(1.3)}), we assume 
$$
F(x,t) = q_1(x,t)f_1(x), \quad G(x,t) = q_2(x,t)f_2(x), \quad (x,t) \in Q,
$$
where \begin{equation}\label{lopukh5}q_1, q_2 \in W^{1,\infty}(0,T;L^{\infty}(\OOO))
\end{equation}
 are given functions.

In addition to (\ref{lopukh1}) and (\ref{1.7}), we assume
\begin{eqnarray}\label{lopukh2}
\ppp_ta_{ij}, \, \ppp_tb_{ij} \in C^1(\ooo{Q}), \quad
\ppp_ta_k, \ppp_tb_k \in L^{\infty}(Q)
\quad \mbox{for $1\le i,j\le d$ and $0\le k \le d$},\nonumber \\
\ppp_t\widetilde a_{\gamma} \in L^{\infty}(Q) \quad 
\mbox{for $\vert \gamma \vert \le 2$}, \quad \ppp_tp, \ppp_tq 
\in C^1(\ppp\OOO\times [0,T]).
\end{eqnarray}
We arbitrarily fix $t_0 \in (0,T)$ and a non-empty open interval 
$I \in (0,T)$ such that $t_0 \in I$.  Furthermore we are given
\begin{equation}\label{lopukh4}
u_0(x) := u(x,t_0), \quad v_0(x):= v(x,t_0), \quad x\in \OOO.
\end{equation}
Then
\\
\begin{thm}\label{t2} {(global unconditional Lipschitz stability for an 
inverse source problem)}
{\it 
Assume that (\ref{lopukh1}), (\ref{lopukh}), (\ref{1.7}), (\ref{lopukh2}) and (\ref{lopukh5}).
Let $\partial_t c_0\in L^\infty(Q), g,h, \partial_t g,\partial_t h\in L^2(I; H^{\hhalf}(\ppp\OOO))\cap H^1(\Gamma \times I)$ and  $u, v \in H^{2,1}(Q)$ satisfy (\ref{(1.3)}), (\ref{(1.4)}), (\ref{lopukh4}), $\ppp_tu, \ppp_tv 
\in H^{2,1}(Q)$. 
We assume
\begin{equation}\label{(1.5)}
\vert q_1(x,t_0)\vert>0, \quad \mbox{and}\quad \vert q_2(x,t_0)\vert > 0 \quad \mbox{for all $x \in \ooo{\OOO}$}.
\end{equation}
Then there exists a constant $C>0$ such that 
\begin{align*}
& \Vert f_1\Vert_{L^2(\OOO)} + \Vert f_2\Vert_{L^2(\OOO)}
\le C\biggl(\Vert u(\cdot,t_0)\Vert_{H^2(\OOO)}
+ \Vert v(\cdot,t_0)\Vert_{H^2(\OOO)}\\
+& \sum_{k=0}^1 (\Vert \ppp_t^ku\Vert_{H^1(\Gamma \times I)}
+ \Vert \ppp_t^kv \Vert_{H^1(\Gamma \times I)})\\
+& \sum_{k=0}^1 (\Vert \ppp_t^kg \Vert
_{H^1(I; L^2(\ppp\OOO\setminus \Gamma))} 
+ \Vert \ppp_t^kg \Vert_{L^2(I; H^{\hhalf}(\ppp\OOO))}
+ \Vert \ppp_t^kh \Vert_{H^1(I; L^2(\ppp\OOO\setminus \Gamma))} 
+ \Vert \ppp_t^kh \Vert_{L^2(I; H^{\hhalf}(\ppp\OOO))}) \biggr).
\end{align*}
}
\end{thm}
In particular, in the case of the homogeneous Robin boundary condition, we have
\begin{cor}
{\it
Under the conditions of Theorem \ref{t2}, we assume that $g=h=0$ 
on $\ppp\OOO\times (0,T)$.
Then there exists a constant $C>0$ such that 
\begin{align*}
& \Vert f_1\Vert_{L^2(\OOO)} + \Vert f_2\Vert_{L^2(\OOO)}
\le C\biggl(\Vert u(\cdot,t_0)\Vert_{H^2(\OOO)}
+ \Vert v(\cdot,t_0)\Vert_{H^2(\OOO)}\\
+& \sum_{k=0}^1 (\Vert \ppp_t^ku\Vert_{H^1(\Gamma \times I)}
+ \Vert \ppp_t^kv \Vert_{H^1(\Gamma \times I)})\biggr).
\end{align*}
}
\end{cor}
We emphasize that in our stability estimate, we do not use neither data  
$u(\cdot,T)$ nor $v(\cdot,0)$ in $\OOO$, 
nor any a priori bounds on $u,v,f_1, f_2$, but we require data $u(\cdot,t_0)$ 
and $v(\cdot,t_0)$ over $\OOO$ at an intermediate time $t_0 \in (0,T)$.
Our stability can be understood unconditional in the sense that 
we do not need to assume any a priori boundedness conditions.
\\
\vspace{0.1cm}

Our key is a classical Carleman estimate for a single parabolic equation
with singular weight function by Imanuvilov \cite{Ima}.
The linearized mean field game equations (\ref{(1.3)}) have two features:
\begin{itemize}
\item
The equation in $u$ is backward and the one in $v$ is forward.
\item
The equation in $v$ contains the second-order spatial derivatives 
$A_0(t)u$ of $u$.
\end{itemize}

The mixed forward and backward equations in (\ref{(1.3)}),
makes the forward problem such as initial boundary 
value problem difficult, but thanks to the symmetry of the time
variable in the weight function, this does not matter for 
Carleman estimates, our main tool.

Systems with coupled principal parts usually cause
difficulty for establishing relevant Carleman estimates.
However, in our case, although the second equation in (\ref{(1.3)}) is 
coupled with the second-order terms $A_0(t)u$ of $u$, the first one in 
(\ref{(1.3)}) is coupled only 
with zeroth order term of $v$, which enables us to execute a typical argument
of absorbing the second-order terms of $u$ by taking large parameters
of the Carleman estimate, as is described in Section 2.  

This article is composed of five sections.  In Section \ref{A2}, we prove
a Carleman estimate (Theorem \ref{t3}) for the linearized mean field game
equations and complete the proof of Theorem \ref{t1}.  
Section \ref{Z3} is devoted to the proof of Theorem \ref{t2}.  
In Section \ref{Z4}, we 
consider the state determination problem for the original nonlinear 
mean value field equations (\ref{(1.1)}).
Section \ref{Z5} is devoted to the proof of Lemma \ref{1}: the key Carleman 
estimate for the case
of the non-homogeneous Robin boundary condition. 

\section{Carleman estimate for the mean field game equations and the proof 
of Theorem \ref{t1}}\label{A2}
 
Let $\Gamma \subset \ppp\OOO$ be an arbitrarily given subboundary.
We choose a subboundary $\Gamma_1 \subset \ppp\OOO$ such that
$\ooo{\ppp\OOO \setminus \Gamma}\subset \Gamma_1$.  Then it is known 
(e.g., Fursikov and Imanuvilov \cite{FI}, Imanuvilov \cite{Ima}) that
there exists a function $\eta \in C^1(\ooo{\Omega})$ such that 
$$
\eta(x)>0\quad\mbox{in}\quad \Omega, \quad \eta\vert_{\Gamma_1} = 0, \quad
\nabla \eta \ne 0\quad \mbox{on}\quad \ooo{\Omega}. 
$$
Let function $\mu = \mu(t)$ satisfy
\begin{equation}\label{(2.1)}
\left\{\begin{array}{rl}
& \mu \in C^{\infty}[0,T], \quad \mu(t) = t^2 \quad \mbox{for $0\le
t \le \frac{T}{4}$}, \\
& \mbox{$\mu(t)$ is monotone increasing in $\left[0,\, \frac{T}{2}\right]$},
\quad \mu(t) = \mu(T-t) \quad \mbox{for $0\le t\le T$}.
\end{array}\right.
\end{equation}
For arbitrarily chosen sufficiently large constant $\la > 0$, we set
$$
\va(x,t) = \frac{e^{\la\eta(x)}}{\mu(t)}, \quad
\alpha(x,t) = \frac{e^{\lambda\eta(x)} 
- e^{2\lambda ||\eta||_{C(\ooo{\Omega})}} }{\mu(t)}, \quad (x,t)\in 
Q:= \OOO \times (0,T).
$$
We recall the assumptions (\ref{lopukh1}) and (\ref{1.7}), and 
further set
\begin{equation}\label{(2.2)}
M:= \sum_{i,j=1}^d ||a_{ij}||_{C^1(\overline{Q})} 
+ \sum_{j=1}^d ||a_j||_{L^{\infty}(Q)}
+ ||c_0||_{L^\infty(Q)}, \quad 
M_0:= \sum_{i,j=1}^d ||a_{ij}||_{C^1(\overline{Q})}.  
\end{equation}

We consider a
boundary value problem
\begin{equation}\label{(2.3)}
\ppp_tu + A(t)u = F \quad \mbox{in $Q$} \quad \mbox{or}
\quad \ppp_tu - B(t)u = F \quad \mbox{in $Q$},  
\end{equation}
and
\begin{equation}\label{(2.4)}
\ppp_{\nu_A}u - p(x,t)u = g \quad \mbox{on $\ppp\OOO\times (0,T)$}.
\end{equation}

Now we state the key Carleman estimate for a parabolic equation.
\\
\begin{lem} \label{1} 
{\it Assume that (\ref{lopukh1}), (\ref{lopukh}) and $p \in C^1(\ppp\OOO\times [0,T])$.
We choose sufficiently large $\la>0$.
Let $g\in L^2(0,T,H^\frac 12 (\partial\Omega))$, 
$\partial_tg\in L^2((\ppp\OOO\setminus \Gamma) \times (0,T)), F\in L^2(Q)$.
Then there exist constants $s_0 > 0$ and $C>0$ independent of $u$ 
such that
\begin{equation}\label{(2.5)}
\int_Q \left(
\frac{1}{{s}\varphi} \left(\vert\partial_t u\vert^2 
+ \sum_{i,j=1}^d\vert \partial_i\ppp_ju\vert^2\right) 
+ s\varphi |\nabla u|^2 + s^3\varphi^3 \vert u\vert^2
\right) e^{2s\alpha}  dx\,dt
\end{equation}
\begin{align*}
\le & C\biggl(
\int_Q \vert F\vert^2 \weight dxdt
+ \int_{(\ppp\OOO\setminus \Gamma)\times (0,T)}
\left(\frac{\vert \partial_t g\vert^2}{s^2\varphi^2}
+ \frac{1}{\root\of{s\varphi}}\vert g\vert^2\right) e^{2{s}\alpha} dSdt
+ \Vert ge^{s\alpha}\Vert^2_{L^2(0,T;H^\frac 12(\partial\Omega))}
\biggr)\\
+ & C\int_{\Gamma \times (0,T)} \left(s\va\vert \nabla u\vert^2
+ s^3\va^3\vert u\vert^2 + \frac{\vert \ppp_tu\vert^2}{s\va}\right)
\weight dSdt 
\end{align*}
for each $s \ge s_0$ and $u\in H^{2,1}(Q)$ 
satisfying (\ref{(2.3)}) and (\ref{(2.4)}).
Here the constant $C>0$ depends continuously on 
$\Vert p\Vert_{C^1(\ppp\OOO\times [0,T])}$,
$M$, $\lambda$. 
}
\end{lem}

Here and henceforth, by $C>0$ we denote generic constants which 
are independent of the parameter $s>0$, and we write $\www{C}(s)$ 
when we need to specify the dependency.

In the case where $g\equiv 0$ on $\ppp\OOO\times (0,T)$ in (\ref{(2.4)}), 
the proof is done
similarly to Fursikov and Imanuvilov \cite{FI}, Imanuvilov \cite{Ima}, which 
treat the zero Dirichlet boundary condition, and Chae, Imanuvilov and Kim
\cite{CIK} for the zero Neumann boundary condition.
For the case $g\not\equiv 0$ with Robin boundary condition, we need to 
modify the proof.  For completeness, we provide the proof of Lemma \ref{1}
in Section \ref{Z5}.
\\

Moreover, we set 
\begin{equation}\label{(2.6)} \Vert g\Vert_*:= \Vert g\Vert_{H^1(0,T;L^2(\ppp\OOO\setminus \Gamma))}
+ \Vert g\Vert_{L^2(0,T;H^{\hhalf}(\ppp\OOO))}.
\end{equation}

We note 
$$
\ppp_i\va = \la(\ppp_i\eta)\va, \quad
\ppp_i\ppp_j\va = (\la\ppp_i\ppp_j\eta + \la^2(\ppp_i\eta)(\ppp_j\eta))\va,
$$
for $1\le i,j \le d$ and
$$
\left\vert \frac{d}{dt}\left(\frac{1}{\mu(t)}\right) \right\vert
\le \frac{C}{\mu^2(t)}, \quad 0<t<T.
$$
Hence 
\begin{equation}\label{(2.10)}
\left\{ \begin{array}{rl}
& \vert \ppp_t\va\vert \le C\va^2, \quad 
\vert \nabla\va\vert \le C\va, \quad
\vert \ppp_i\ppp_j\va\vert \le C\va \quad \mbox{in $Q$ for $1\le i,j\le d$},
\cr\\
& \vert \nabla \alpha \vert \le C\va, \quad \vert \ppp_t\alpha\vert 
\le C\va^2 \quad \mbox{in $Q$} \cr
\end{array}\right.
\end{equation}

In order to rewrite the norms appearing in (\ref{(2.5)}), we show
\\
\begin{lem}\label{3}
{\it
(i) For each $\rho\in \R$, we have
\begin{equation}\label{(2.7)}
\sup_{s\ge 1}
\sup_{(x,t)\in Q} \vert \va(x)^{\rho}e^{2s\alpha(x,t)}\vert < \infty.
\end{equation}
\\
(ii) Let $\rho\in \R$ and $\psi \in C([0,T]; C^1(\ppp\OOO))$
be arbitrarily given function and let $s\ge 1$ be arbitrary.  Then there exist a constant $C_{\rho},$ independent of $s$, such that
\begin{equation}\label{(2.8)}
\Vert \va^{\rho}\psi ge^{s\alpha}\Vert_* \le C_{\rho}s\Vert g\Vert_*.
\end{equation}
}
\end{lem}

{\bf Proof of Lemma \ref{3}.}
\\
(i) First we have
$$
\va(x,t)^{\rho}e^{2s\alpha} \le \va(x,t)^{\rho}
\le \mu(t)^{-\rho} e^{\la\rho\eta(x)},
\quad (x,t) \in Q
$$
if $\rho\le 0$, which readily verifies (\ref{(2.7)}).
On the other hand, for $\rho>0$, we have
$$
\va(x,t)^{\rho}e^{2s\alpha} 
\le \frac{C}{\mu(t)^{\rho}} 
\exp\left( 
2 \frac{e^{\la\eta(x)} - e^{2\la\Vert\eta\Vert_{C(\ooo{\OOO})}} }{\mu(t)} 
\right)
\le \frac{C}{\mu(t)^{\rho}}e^{-\frac{C_1}{\mu(t)}}, \quad (x,t) \in Q,
$$
where 
$$
C_1:=  2(e^{2\la\Vert\eta\Vert_{C(\ooo{\OOO})}}
- e^{\la\Vert\eta\Vert_{C(\ooo{\OOO})}}) > 0.
$$
Noting that $\xi:= \frac{1}{\mu(t)}$, $t\in (0,T)$ varies $[c_0,\infty)$ 
with some constant
$c_0>0$ and $\sup_{\xi\ge c_0} \xi^{\rho}e^{-C_1\xi} < \infty$, 
we see (\ref{(2.7)}) for all $\rho \in \R$.
\\
(ii) By (\ref{(2.10)}), we have
\begin{align*}
& \Vert \va^{\rho}\psi ge^{s\alpha}\Vert_{H^1(0,T;L^2(\ppp\OOO\setminus
\Gamma))}
\le \Vert \rho\va^{\rho-1}(\ppp_t\va) \psi ge^{s\alpha}\Vert
_{L^2(0,T;L^2(\ppp\OOO\setminus \Gamma))}\\
+& \Vert \va^{\rho}(\ppp_t\psi)ge^{s\alpha}\Vert
_{L^2(0,T;L^2(\ppp\OOO\setminus\Gamma))}
+ \Vert \va^{\rho}\psi (\ppp_tg)e^{s\alpha}\Vert
_{L^2(0,T;L^2(\ppp\OOO\setminus \Gamma))}\\
+& s\Vert \va^{\rho}\psi g(\ppp_t\alpha)e^{s\alpha}\Vert
_{L^2(0,T;L^2(\ppp\OOO\setminus\Gamma))}
+ \Vert \va^{\rho}\psi ge^{s\alpha}\Vert_{L^2(0,T;L^2(\ppp\OOO\setminus
\Gamma))}\\
\le& Cs\Vert \va^{\rho+2}ge^{s\alpha}\Vert_{L^2(0,T;L^2(\ppp\OOO \setminus
\Gamma))} 
+ C\Vert \va^{\rho}(\ppp_tg)e^{s\alpha}\Vert
_{L^2(0,T;L^2(\ppp\OOO \setminus\Gamma))}
+ C\Vert \va^{\rho+1}ge^{s\alpha}\Vert_{L^2(0,T;L^2(\ppp\OOO \setminus
\Gamma))}.
\end{align*}
Therefore, (\ref{(2.7)}) yields
\begin{align*}
& \Vert \va^{\rho}\psi ge^{s\alpha}\Vert_{H^1(0,T;L^2(\ppp\OOO
\setminus \Gamma))}\\
\le& C(s\Vert g\Vert_{L^2(0,T;L^2(\ppp\OOO\setminus \Gamma))}
+ \Vert \ppp_tg\Vert_{L^2(0,T;L^2(\ppp\OOO\setminus \Gamma))})
\le Cs\Vert g\Vert_*
\end{align*}
for all $s\ge 1$.
 
Next, in view of the Sobolev-Slobodecki norm in $H^{\hhalf}(\ppp\OOO)$ 
(e.g., Adams \cite{Ad}), we can directly verify  that there exist a constant $C$ such that
\begin{equation}\label{(2.9)}
\Vert \psi a\Vert_{H^{\hhalf}(\ppp\OOO)}
\le C\Vert \psi\Vert_{C^1(\ppp\OOO)} \Vert a\Vert_{H^{\hhalf}(\ppp\OOO)}
                                           \end{equation}
for $a \in H^{\hhalf}(\ppp\OOO)$ and $\psi \in C^1(\ppp\OOO)$.

By (\ref{(2.7)}), we have
$$
\Vert \va^{\rho}\psi e^{s\alpha}(\cdot,t)\Vert_{C^1(\ooo{\OOO})}
\le \Vert \va^{\rho}\psi e^{s\alpha}(\cdot,t)\Vert_{C(\ooo{\OOO})}
+ \Vert \nabla(\va^{\rho}\psi e^{s\alpha})(\cdot,t)\Vert_{C^1(\ooo{\OOO})}
\le C_{\rho}s.
$$
Hence, (\ref{(2.10)}) and (\ref{(2.9)}) yield
$$
\Vert \va^{\rho}\psi g e^{s\alpha}(\cdot,t)\Vert_{H^{\hhalf}(\ppp\OOO)}
\le C\Vert \va^{\rho}\psi e^{s\alpha}(\cdot,t)\Vert_{C^1(\ooo{\ppp\OOO})}
\Vert g(\cdot,t)\Vert_{H^{\hhalf}(\ppp\OOO)}
\le C_{\rho}s\Vert g(\cdot,t)\Vert_{H^{\hhalf}(\ppp\OOO)},
$$
and so 
$$
\Vert \va^{\rho}\psi g e^{s\alpha}\Vert_{L^2(0,T;H^{\hhalf}(\ppp\OOO))}
\le C_{\rho}s\Vert g \Vert_{L^2(0,T;H^{\hhalf}(\ppp\OOO))}.
$$
Thus the proof of (\ref{(2.8)}) is complete.
$\blacksquare$

In view of Lemma \ref{3}, we can rewrite (\ref{(2.5)}) as 
$$
\int_Q \left(
\frac{1}{s\varphi} \left(\vert\partial_t u\vert^2 
+ \sum_{i,j=1}^d\vert \partial_i\ppp_ju\vert^2\right) 
+ s\varphi |\nabla u|^2 + s^3\varphi^3 \vert u\vert^2
\right) e^{2s\alpha}  dx\,dt
$$
\begin{equation}\label{2.10}
\le C\int_Q \vert F\vert^2 \weight dxdt
+ \www{C}(s)(\Vert g\Vert^2_* + \Vert u\Vert_{H^1(\Gamma \times (0,T))}^2)
\end{equation}
for each $s \ge s_0$.
\\

Next we prove
\begin{lem}\label{2}
{\it 
Assume (\ref{lopukh1}), (\ref{lopukh}) and $p\in C^1(\partial\Omega\times [0,T]).$
Let  $g\in L^2(0,T,H^\frac 12 (\partial\Omega))$, 
$\partial_tg\in L^2((\ppp\OOO\setminus \Gamma) \times (0,T)), F\in L^2(Q)$.  
We fix a sufficiently large constant $\la>0$.
Then, for each $m\in \R$,  there exist constants $s_0 > 0$ and $C>0$ such that 
\begin{align*}
& \int_Q \left((s\va)^{m-1}\left(\vert \ppp_tu\vert^2 
+ \sumij \vert \ppp_i\ppp_ju\vert^2\right)
+ (s\va)^{m+1}\vert \nabla u\vert^2 + (s\va)^{m+3}\vert u\vert^2\right)
\weight dxdt\\
\le& C\int_Q (s\va)^m\vert F\vert^2 \weight dxdt\\
+ & C\int_{(\ppp\OOO\setminus \Gamma)\times (0,T)}
s^m\left( \frac{\vert \ppp_t(\va^{\frac{m}{2}}g)\vert^2}{s^2\va^2} 
+\frac{1}{\root\of{s\varphi}} \vert \va^{\frac{m}{2}} g\vert^2 \right) e^{2s\alpha} dSdt\\
+ & Cs^m \Vert \va^{\frac{m}{2}} ge^{s\alpha}\Vert^2_{L^2(0,T;
H^{\hhalf}(\ppp\OOO))} + \www{C}(s,m)\Vert u\Vert_{H^1(\Gamma \times (0,T))}^2
\end{align*}
for all $s > s_0$ and $u\in H^{2,1}(Q)$ satisfying (\ref{(2.3)}) and 
(\ref{(2.4)}).
}
\end{lem}

Here $\www{C}(s,m)$ is a positive constant depending on 
$s$ and $m$.
\\
{\bf Proof of Lemma \ref{2}.}
\\
We will derive Lemma \ref{2} from Lemma \ref{1} in the case where 
$\ppp_tu + A(t)u = F$.  The derivation for the case of
$\ppp_tu - B(t)u = F$ is quite similar.

Moreover, choosing $s>0$ sufficiently large, by 
$a_j, b_j, a_0, b_0 \in L^{\infty}(Q)$, we can absorb the lower-order terms 
$\sum_{j=1}^d a_j\ppp_ju$, $a_0u$, $\sum_{j=1}^d b_j\ppp_ju$, $b_0u$
into the left-hand side of the Carleman estimate.
Thus it suffices to prove the lemma for 
$A(t) = \sumij a_{ij}\ppp_i\ppp_j$ and $B(t) = \sumij b_{ij}\ppp_i\ppp_j$.
\\

We set 
$$
w:= \va^{\frac{m}{2}}u.
$$
Then we can directly calculate:
\begin{equation}\label{(2.11)}
\left\{ \begin{array}{rl}
& \ppp_tw = \frac{m}{2}\va^{\frac{m}{2}-1}(\ppp_t\va)u
+ \va^{\frac{m}{2}}\ppp_tu,\\
& \ppp_iw = \frac{m}{2}\va^{\frac{m}{2}-1}(\ppp_i\va)u
+ \va^{\frac{m}{2}}\ppp_iu, \\
& \ppp_i\ppp_jw = \frac{m}{2}\left(\frac{m}{2}-1\right)\va^{\frac{m}{2}-2}
(\ppp_i\va)(\ppp_j\va)u
+ \frac{m}{2}\va^{\frac{m}{2}-1}(\ppp_i\ppp_j\va)u\\
+ &\frac{m}{2}\va^{\frac{m}{2}-1}(\ppp_i\va)(\ppp_ju)
+ \frac{m}{2}\va^{\frac{m}{2}-1}(\ppp_j\va)(\ppp_iu)
+ \va^{\frac{m}{2}}\ppp_i\ppp_ju.
\end{array}\right.
                                             \end{equation}
Therefore,
\begin{align*}
& \ppp_tw + A(t)w 
= \va^{\frac{m}{2}}F + \frac{m}{2}\va^{\frac{m}{2}-1}(\ppp_t\va)u\\
+& \frac{m}{2}\left( \frac{m}{2}-1 \right)\sumij a_{ij}
  \va^{\frac{m}{2}-2}(\ppp_i\va)(\ppp_j\va)u
+ \frac{m}{2} \va^{\frac{m}{2}-1}\sumij a_{ij}(\ppp_i\ppp_j\va)u\\
+ &m\va^{\frac{m}{2}-1}\sumij a_{ij}(\ppp_i\va)(\ppp_ju)  \quad \mbox{in $Q$}
\end{align*}
and
\begin{align*}
& \ppp_{\nu_A}w = \sumij a_{ij}\ppp_i(\va^{\frac{m}{2}}u)\nu_j
= \va^{\frac{m}{2}}\sumij a_{ij} (\ppp_iu)\nu_j 
+ \frac{m}{2} \va^{\frac{m}{2}-1}\left(\sumij a_{ij}(\ppp_i\va)\nu_j
\right)u\\
=& \va^{\frac{m}{2}}g 
+ \left( \frac{m}{2}\va^{-1}\sumij a_{ij} (\ppp_i\va)\nu_j \right)w
\quad \mbox{on $\ppp\OOO \times (0,T)$}.
\end{align*}
Hence,
\begin{equation}\label{1313}
\ppp_{\nu_A}w - \www{p}(x,t)w = \va^{\frac{m}{2}}g, \quad
(x,t) \in \ppp\OOO\times (0,T),
\end{equation}
where 
$$
\www{p}(x,t):= p(x,t) + \frac{m}{2}\va^{-1}\sumij a_{ij} (\ppp_i\va)\nu_j
= p(x,t) + \frac{m\la}{2}\sumij a_{ij}(\ppp_i\eta)\nu_j.
$$
Moreover 
\begin{equation}\label{1414}
\ppp_tw + A(t)w = \va^{\frac{m}{2}}F + \www{F} \quad \mbox{in $Q$},
\end{equation}
where we see
$$
\vert \www{F}(x,t)\vert \le C(\va^{\frac{m}{2}+1}\vert u\vert 
+ \va^{\frac{m}{2}}\vert \nabla u(x,t)\vert 
\le C(\va\vert w(x,t)\vert + \vert \nabla w(x,t)\vert), \quad
(x,t) \in Q.
$$
We apply Lemma \ref{1} to (\ref{1313}) and (\ref{1414}) to 
obtain
\begin{eqnarray}\label{(2.12)}
 \int_Q \left( \frac{1}{s\va}\left( \vert \ppp_tw\vert^2 
+ \sumij \vert \ppp_i\ppp_jw\vert^2\right)
+ s\va\vert \nabla w\vert^2 + s^3\va^3\vert w\vert^2\right)
\weight dxdt\nonumber\\
\le C\int_Q \va^m\vert F\vert^2 \weight dxdt
+ C\int_Q (\va^2\vert w\vert^2 + \vert \nabla w\vert^2) \weight dxdt \nonumber\\
+  C\int_{(\ppp\OOO\setminus \Gamma)\times (0,T)}
\left( \frac{\vert \ppp_t(\va^{\frac{m}{2}}g)\vert^2}{s^2\va^2} 
+ \frac{1}{\sqrt{s\va}}\vert \va^{\frac{m}{2}}g\vert^2 \right) 
e^{2s\alpha} dSdt\nonumber\\
+  C\Vert \va^{\frac{m}{2}} ge^{s\alpha}\Vert^2_{L^2(0,T;
H^{\hhalf}(\ppp\OOO))} + \www{C}(s)\Vert w\Vert_{H^1(\Gamma \times (0,T))}^2
\end{eqnarray}
for all $s > s_0$.
Choosing $s>0$ large, we can absorb the second term on the right-hand side
into the left-hand side.  

In terms of $u$ we rewrite as follows.
By (\ref{(2.11)}), we first have
$$
s^3\va^3\vert w\vert^2 = s^3\va^{m+3}\vert u\vert^2,
$$
$$
\va^{\frac{m}{2}}(\ppp_iu) = \ppp_iw - \frac{m}{2}\va^{\frac{m}{2}-1}
(\ppp_i\va)u, 
$$
and so
$$
s\va^{m+1}\vert \nabla u\vert^2 \le 2s\va\vert \nabla w\vert^2
+ Cs\va\left\vert m\va^{\frac{m}{2}-1}(\nabla \va)u
\right\vert^2
\le 2s\va\vert \nabla w\vert^2 + Cs^2\va^{m+2}\vert u\vert^2.
$$
Moreover, (\ref{(2.11)}) implies
$$
\frac{1}{s} \va^{m-1}\vert \ppp_tu\vert^2
\le \frac{C}{s\va}\vert \ppp_tw\vert^2 + Cs^{-1}\va^{m+1} \vert u\vert^2.
$$
Finally, again (\ref{(2.11)}) yields
\begin{align*}
& \va^{\frac{m}{2}}\ppp_i\ppp_ju  
= \ppp_i\ppp_jw - \frac{m}{2}\left( \frac{m}{2}-1\right) \va^{\frac{m}{2}-2}
(\ppp_i\va)(\ppp_j\va)u
- \frac{m}{2}\va^{\frac{m}{2}-1}(\ppp_i\ppp_j\va)u\\
-& \frac{m}{2}\va^{\frac{m}{2}-1}(\ppp_i\va)(\ppp_ju)
- \frac{m}{2}\va^{\frac{m}{2}-1}(\ppp_j\va)(\ppp_iu),
\end{align*}
and so
\begin{align*}
& \frac{1}{s\va} \left\vert \va^{\frac{m}{2}}\ppp_i\ppp_ju\right\vert^2
\le \frac{C}{s\va}\vert \ppp_i\ppp_jw\vert^2\\
+ & \frac{C}{s\va}\left\vert \frac{m}{2}\left( \frac{m}{2} - 1\right)
\va^{\frac{m}{2}-2}(\ppp_i\va)(\ppp_j\va)u
+ \frac{m}{2}\va^{\frac{m}{2}-1}(\ppp_i\ppp_j\va)u
+ \frac{m}{2} \va^{\frac{m}{2}-1}((\ppp_i\va)\ppp_ju
  + (\ppp_j\va)\ppp_iu)\right\vert^2.
\end{align*}
Hence,
$$
\frac{1}{s} \va^{m-1}\vert \ppp_i\ppp_ju\vert^2
\le \frac{C}{s\va}\vert \ppp_i\ppp_jw\vert^2
+ \frac{C}{s\va}(\va^m\vert u\vert^2 + \va^m\vert \nabla u\vert^2).
$$
In (\ref{(2.12)}), we can estimate $\Vert w\Vert_{H^1(\Gamma \times (0,T))}^2$ by means of Lemma \ref{3} (i):
$$
 \Vert w\Vert_{H^1(\Gamma \times (0,T))}^2\le \www{C}(s,m)\Vert u\Vert_{H^1(\Gamma \times (0,T))}^2
$$
for all $s\ge 1$.  Thus, the proof of Lemma \ref{2} is complete 
$\blacksquare$
\\

In particular, setting $m=1$ in Lemma \ref{2}, we have
\\
\begin{lem}\label{4}
{\it 
Suppose that conditions of Lemma \ref{2} holds true.
We can find constants $s_0 > 0$ and $C>0$ such that 
\begin{align*}
& \int_Q \left(\vert \ppp_tu\vert^2 + \sumij \vert \ppp_i\ppp_ju\vert^2
+ s^2\va^2\vert \nabla u\vert^2 + s^4\va^4\vert u\vert^2\right)
\weight dxdt\\
\le& C\int_Q s\va\vert \ppp_tu + A(t)u\vert^2 \weight dxdt\\
+ & C\int_{(\ppp\OOO\setminus \Gamma)\times (0,T)}
\left( \frac{\vert \ppp_t(\va^{\hhalf}g)\vert^2}{s^2\va^2} 
+ {\root\of{s\varphi}}\vert g\vert^2 \right) e^{2s\alpha} dSdt\\
+ & C\Vert s^{\hhalf}\va^{\hhalf} ge^{s\alpha}\Vert^2_{L^2(0,T;
H^{\hhalf}(\ppp\OOO))} + \www{C}(s)+\Vert u\Vert_{H^1(\Gamma \times (0,T))}^2
\end{align*}
for all $s > s_0$ and $u\in H^{2,1}(Q)$ satisfying (\ref{(2.3)}),
(\ref{(2.4)}).
}\end{lem}

Now, noting that we have the Carleman estimates both for 
$\ppp_t + A(t)$ and $\ppp_t - B(t)$ with the same weight $\weight$,
we derive the main Carleman estimate for the mean field game
system (\ref{(1.3)}) with (\ref{(1.4)}).
Setting $F:= c_0v+F$ in the first equation in (\ref{(1.4)}), we apply Lemma \ref{4} to obtain
\begin{equation}\label{(2.13)}
 \int_Q \left(\vert \ppp_tu\vert^2 + \sumij \vert \ppp_i\ppp_ju\vert^2
+ s^2\va^2\vert \nabla u\vert^2 + s^4\va^4\vert u\vert^2 \right)
\weight dxdt
                                     \end{equation}
\begin{eqnarray}
\le  C\int_Q s\va\vert v\vert^2 \weight dxdt 
+ C\int_Q s\va \vert F\vert^2 \weight dxdt\nonumber\\
+  C\int_{(\ppp\OOO\setminus \Gamma)\times (0,T)}
\left( \frac{\vert \ppp_tg\vert^2}{s\va} + \root\of{s\va}\vert g\vert^2
\right) \weight dSdt
+ s\Vert \va^{\frac{1}{2}} ge^{s\alpha}\Vert^2_{L^2(0,T;
H^{\hhalf}(\ppp\OOO))} + \www{C}(s) \Vert u\Vert_{H^1(\Gamma \times (0,T))}^2       \nonumber\\
\le C\int_Q s\va\vert v\vert^2 \weight dxdt 
+ C\int_Q s\va \vert F\vert^2 \weight dxdt
+ \www{C}(s)(\Vert g\Vert_*^2 + \Vert u\Vert_{H^1(\Gamma \times (0,T))}^2)\nonumber
\end{eqnarray}
for all $s > s_0$.

The application of Lemma \ref{1} to the second equation in (\ref{(1.5)}) with 
$G:= G + A_0(t)u$ yields
\begin{eqnarray}\label{(2.14)}
 \int_Q \biggl\{ \frac{1}{s\va}\left(
\vert \ppp_tv\vert^2 + \sumij \vert \ppp_i\ppp_jv\vert^2\right)
+ s\va\vert \nabla v\vert^2 + s^3\va^3\vert v\vert^2 \biggr\}
\weight dxdt\\
\le C\int_Q \sumij \vert \ppp_i\ppp_ju\vert^2 \weight dxdt
+ C\int_Q \vert G\vert^2 \weight dxdt 
+ \www{C}(s)(\Vert h\Vert_*^2 + \Vert v\Vert_{H^1(\Gamma \times (0,T))}^2).\nonumber
                                                \end{eqnarray}
Estimating the first term on the right-hand side of (\ref{(2.14)}) in 
terms of (\ref{(2.13)}), we obtain
\begin{eqnarray}
 \int_Q \biggl\{ \frac{1}{s\va}\left(
\vert \ppp_tv\vert^2 + \sumij \vert \ppp_i\ppp_jv\vert^2\right)
+ s\va\vert \nabla v\vert^2 + s^3\va^3\vert v\vert^2 \biggr\}
\weight dxdt\nonumber\\
\le C\int_Q s\va\vert v\vert^2 \weight dxdt
+ C\int_Q s\va\vert F\vert^2 \weight dxdt 
+ C\int_Q \vert G\vert^2 \weight dxdt \nonumber\\
+  \www{C}(s)(\Vert g\Vert_*^2 + \Vert h\Vert^2_* + \Vert u\Vert_{H^1(\Gamma \times (0,T))}^2 + \Vert v\Vert_{H^1(\Gamma \times (0,T))}^2)\nonumber
\end{eqnarray}
for all large $s > 0$.
Hence, choosing $s>0$ sufficiently large, we can absorb 
the first term on the right-hand side into the left-hand side, and we can 
obtain
\begin{align*}
& \int_Q \biggl\{ \frac{1}{s\va}\left(
\vert \ppp_tv\vert^2 + \sumij \vert \ppp_i\ppp_jv\vert^2\right)
+ s\va\vert \nabla v\vert^2 + s^3\va^3\vert v\vert^2 \biggr\}
\weight dxdt\\
\le & C\int_Q (s\va\vert F\vert^2 + \vert G\vert^2) \weight dxdt
+ \www{C}(s)(\Vert g\Vert_*^2 + \Vert h\Vert^2_* +  \Vert u\Vert_{H^1(\Gamma \times (0,T))}^2 + \Vert v\Vert_{H^1(\Gamma \times (0,T))}^2)
\end{align*}
for all large $s > 0$.  Adding with (\ref{(2.13)}), we absorb the term
$\int_Q s\va \vert v\vert^2 \weight dxdt$ on the right-hand side into the 
left-hand side, so that we proved
\\
\begin{thm}\label{t3} (Carleman estimate for a generalized mean field game
equations)
{\it 
Let $ g,h\in L^2(0,T;H^\frac 12(\partial\Omega)), \partial_t g,\partial_t h\in 
L^2(0,T;L^2(\partial\Omega\setminus \Gamma)), F,G\in L^2(Q)$ and (\ref{lopukh1}),
(\ref{lopukh}),  (\ref{1.7}) holds true.
We fix $\la>0$ sufficiently large.  
Then we can find constants $s_0 > 0$ and $C>0$ such that 
\begin{eqnarray}\label{(2.17)}
 \int_Q \biggl\{ \vert \ppp_tu\vert^2 + \sumij \vert \ppp_i\ppp_ju\vert^2
+ s^2\va^2\vert \nabla u\vert^2 + s^4\va^4\vert u\vert^2 \\
+ \frac{1}{s\va}\left(
\vert \ppp_tv\vert^2 + \sumij \vert \ppp_i\ppp_jv\vert^2\right)
+ s\va\vert \nabla v\vert^2 + s^3\va^3\vert v\vert^2 \biggr\}
\weight dxdt\nonumber\\
\le C\int_Q (s\va\vert F\vert^2 + \vert G\vert^2) \weight dxdt
+ \www{C}(s)(\Vert g\Vert_*^2 + \Vert h\Vert^2_* +  \Vert u\Vert_{H^1(\Gamma \times (0,T))}^2 + \Vert v\Vert_{H^1(\Gamma \times (0,T))}^2)\nonumber
\end{eqnarray}
for all $s > s_0$ and $u, v\in H^{2,1}(Q)$ satisfying
(\ref{(1.3)}) and (\ref{(1.4)}).
Here the constant $C>0$ depends continuously on $M$: bound of the
coefficients and $\la$ but
independent of $s \ge s_0$.
}\end{thm}

Now we proceed to 
\\
{\bf Proof of Theorem \ref{t1}.}
\\
Since, $\mu(t) \ge \mu(\ep)$ for $\ep \le t \le T-\ep$, for some positive constants $C_2$ and $C_3$  we have
$$
\alpha(x,t) \ge \frac{e^{\la\eta(x)} - e^{2\la\Vert \eta\Vert_{C(\ooo{\OOO})}} 
}{\mu(\ep)}
\ge \frac{-C_2}{\mu(\ep)} =: -C_3 
$$
for all $x\in \ooo{\OOO}$ and $\ep \le t \le T-\ep$, we see that 
$$
e^{2s\alpha(x,t)} \ge e^{-2sC_3}, \quad x\in \ooo{\OOO}, \, 
\ep\le t \le T-\ep.
$$
Thus Theorem \ref{t3} completes the proof of Theorem \ref{t1}.
$\blacksquare$
\section{Proof of Theorem \ref{t2}}\label{Z3}

The proof is based on a similar idea to Theorem 3.1 
in Imanuvilov and Yamamoto \cite{IY98}, where we have to 
estimate extra second-order derivatives of $u$. 

Without loss of generality, we can assume that $t_0 = \frac{T}{2}$ by 
scaling the time variable.
Indeed, we choose small $\delta > 0$ 
such that $0<t_0-\delta < t_0 < t_0 + \delta < T$.
Then we consider a change of the variables 
$t \mapsto \xi:= \frac{t-(t_0-\delta)}{2\delta}T$.  Then, the inverse problem 
over the time interval $(t_0-\delta,\, t_0+\delta)$ can be transformed to 
$(0,T)$ with $t_0= \frac{T}{2}$.  
Thus, it is sufficient to assume that $t_0 = \frac{T}{2}$ and
$I = (0,T)$.
\\
{\bf First Step.}
\\
We show
\\
\begin{lem}\label{5}
{\it 
Let $r \in \R$ and $w\in H^2(Q)$.  Then there exist a constant $C>0$  independent of $w$ such that 
\begin{equation}\label{(3.1)}
\int_{\ppp\OOO} \va^{2r}\vert w\vert^2 e^{2s\alpha} dS
\le C\int_{\OOO}
(\va^{2r}\vert \nabla w\vert^2 + s^2\va^{2r+2}\vert w\vert^2)\weight dx
                                     \end{equation}
and
\begin{equation}\label{(3.2)}
\int_{\ppp\OOO} \va^{2r}\vert \nabla w\vert^2 e^{2s\alpha} dS
\le C\int_{\OOO} \left(\va^{2r}\sumij\vert \ppp_i\ppp_j w\vert^2 
+ s^2\va^{2r+2}\vert \nabla w\vert^2\right)\weight dx             \end{equation}
for all $s>0$.
}
\end{lem}
{\bf Proof of Lemma \ref{5}.}
\\
Indeed, the trace theorem and (\ref{(2.10)}) imply
\begin{align*}
& \Vert \va^rw e^{s\alpha}\Vert^2_{L^2(\ppp\OOO)}
\le C(\Vert \va^r we^{s\alpha}\Vert^2_{L^2(\OOO)}
+ \Vert \nabla(\va^r we^{s\alpha})\Vert^2_{L^2(\OOO)})\\
\le& C\int_{\OOO}
(\va^{2r}\vert w\vert^2 + s^2\va^{2r+2}\vert w\vert^2
+ \va^{2r}\vert \nabla w\vert^2) \weight dx\\
\le& C\int_{\OOO}
(\va^{2r}\vert \nabla w\vert^2 + s^2\va^{2r+2}\vert w\vert^2)\weight dx.
\end{align*}
Thus (\ref{(3.1)}) is seen.  Similarly we can prove (\ref{(3.2)}).
$\blacksquare$
\\
\vspace{0.2cm}
{\bf Second Step: Carleman estimate for $(\ppp_tu, \ppp_tv)$.}
\\
Setting $y:= \ppp_tu$ and $z:= \ppp_t v$ in (\ref{(1.3)}), we have
\begin{equation}\label{(3.3)}
\left\{ \begin{array}{rl}
& \ppp_ty + A(t)y = c_0z + (\ppp_tc_0)v - (\ppp_tA(t))u + \ppp_tF,\\
& \ppp_tz - B(t)z = A_0(t)y + (\ppp_tA_0(t))u + (\ppp_tB(t))v + \ppp_tG
\quad \mbox{in $Q$}
\end{array}\right.
                                      \end{equation}
and
\begin{equation}\label{(3.4)}
\ppp_{\nu_A}y - py = \ppp_tg + g_1, \quad
\ppp_{\nu_B}z - qz = \ppp_th + h_1 \quad \mbox{on $\ppp\OOO \times (0,T)$}.
                                                \end{equation}
Here 
\begin{align*}
& g_1:= \sumij (\ppp_ta_{ij})(\ppp_ju)\nu_i - (\ppp_tp)u, \\
& h_1:= \sumij (\ppp_tb_{ij})(\ppp_jv)\nu_i - (\ppp_tq)v \quad 
\mbox{on $\ppp\OOO \times (0,T)$}.
\end{align*}

In this step, we will prove
\\
\begin{prop}\label{p1} {(Carleman estimate for $(\ppp_tu, \ppp_tv)$)}
{\it Let all the assumptions of Theorem \ref{t2}, except (1.8) and (1.11),  hold true.
There exist constants $s_0>0$ and $C>0$ such that 
\begin{align*}
& \int_Q \biggl\{ \frac{1}{s\va}\left( \vert \ppp_t^2u\vert^2 
+ \sumij \vert \ppp_i\ppp_j\ppp_tu\vert^2\right)
+ s\va\vert \nabla \ppp_tu\vert^2 + s^3\va^3\vert \ppp_tu\vert^2 \\
+& \frac{1}{s^2\va^2}\left( \vert \ppp_t^2v\vert^2 
+ \sumij \vert \ppp_i\ppp_j\ppp_tv\vert^2\right)
+ \vert \nabla \ppp_tv\vert^2 + s^2\va^2\vert \ppp_tv\vert^2
\biggr\} \weight dxdt  \\
\le& C\int_Q \left(s\va\vert F\vert^2 + \vert \ppp_tF\vert^2
+ \vert G\vert^2 + \frac{1}{s\va}\vert \ppp_tG\vert^2\right) 
\weight dxdt\\
+& \www{C}(s)( \Vert \ppp_tg\Vert^2_* + \Vert \ppp_th\Vert^2_*
+ \Vert g\Vert^2_* + \Vert h\Vert^2_*\nonumber\\
+ & \Vert u\Vert_{H^1(\Gamma \times (0,T))}^2 + \Vert v\Vert_{H^1(\Gamma \times (0,T))}^2+ \Vert \partial_tu\Vert_{H^1(\Gamma \times (0,T))}^2 + \Vert \partial_tv\Vert_{H^1(\Gamma \times (0,T))}^2)
\end{align*}
for all large $s > s_0$ and $u,v \in H^{2,1}(Q)$ satisfying 
$\ppp_tu, \ppp_tv \in H^{2,1}(Q)$, (\ref{(1.3)}) and (\ref{(1.4)}).
}
\end{prop}
{\bf Proof of Proposition \ref{p1}.}
\\
Since $A(t)$ and $B(t)$ depend on $t$, after taking the time derivatives of
$u, v$, the first derivatives of $u,v$ enter the Robin boundary conditions 
of $\ppp_tu, \ppp_tv$, and the estimation of $\nabla u, \nabla v$ outside 
the observation subboundary $\Gamma$, is indispensable.
Such estimation can be done also by the Carleman estimate.
 
By (\ref{lopukh2}), we note that $\ppp_tc_0$ and all the coefficients of
$\ppp_tA(t)$, $\ppp_tA_0(t)$ and $\ppp_tB(t)$ are in $L^{\infty}(Q)$.
Therefore, we can apply Lemmata \ref{1} and \ref{2} with $m=-1$ to 
the first and the second equations in (\ref{(3.3)}) respectively, we have
\begin{eqnarray}\label{(3.5)}
 \int_Q \left(\frac{1}{s\va}\left( \vert \ppp_ty\vert^2 
+ \sumij \vert \ppp_i\ppp_jy\vert^2\right)
+ s\va\vert \nabla y\vert^2 + s^3\va^3\vert y\vert^2\right)
\weight dxdt\nonumber\\
\le C\int_Q \vert z\vert^2 \weight dxdt
+ C\int_Q \left( \sum_{\vert \gamma\vert\le 2} \vert \ppp_x^{\gamma}u\vert^2
+ \vert v\vert^2\right) \weight dxdt    \nonumber\\
+ C\int_Q \vert \ppp_tF\vert^2 \weight dxdt + \www{C}(s)\Vert \ppp_tg\Vert^2_*
\\
+ C\int_{(\ppp\OOO\setminus \Gamma)\times (0,T)}
\left( \frac{\vert \ppp_tg_1\vert^2}{s^2\va^2} + \frac{\vert g_1\vert^2}{\root\of{s\varphi}}
\right) e^{2s\alpha} dSdt
+ C\Vert g_1 e^{s\alpha}\Vert^2_{L^2(0,T;H^{\hhalf}(\ppp\OOO))} 
+ \www{C}(s) \Vert y\Vert_{H^1(\Gamma \times (0,T))}^2 \nonumber
\end{eqnarray}
and
\begin{eqnarray}\label{(3.6)}
 \int_Q \left(\frac{1}{s^2\va^2}\left( \vert \ppp_tz\vert^2 
+ \sumij \vert \ppp_i\ppp_jz\vert^2\right)
+ \vert \nabla z\vert^2 + s^2\va^2\vert z\vert^2\right)
\weight dxdt\nonumber\\
\le C\int_Q \frac{1}{s\va}\sumij \vert \ppp_i\ppp_jy\vert^2 \weight dxdt
+ C\int_Q \frac{1}{s\va} \sum_{\vert \gamma\vert\le 2} 
(\vert \ppp_x^{\gamma}u\vert^2 + \vert \ppp_x^{\gamma} v\vert^2) \weight dxdt
                                        \nonumber\\
+  C\int_Q \frac{1}{s\va} \vert \ppp_tG\vert^2 \weight dxdt
+ \www{C}(s)\Vert \va^{-\hhalf} \ppp_th\Vert^2_*            \nonumber \\
+  C\int_{(\ppp\OOO\setminus \Gamma)\times (0,T)}
\left( \frac{\vert \ppp_t(\va^{-\hhalf}h_1)\vert^2}{s^3\va^2} 
+ \frac{s^{-1}\vert \va^{-\hhalf} h_1\vert^2}{\root\of{s\varphi}} \right) e^{2s\alpha} dSdt
\nonumber\\
+ Cs^{-1}\Vert \va^{-\hhalf}h_1 e^{s\alpha}\Vert^2
_{L^2(0,T;H^{\hhalf}(\ppp\OOO))} 
+ \www{C}(s)\Vert z\Vert_{H^1(\Gamma \times (0,T))}^2 .
                             \end{eqnarray}

Applying (\ref{(3.5)}) for estimating the first term on the right-hand side of 
(\ref{(3.6)}), we have
\begin{align*}
& \int_Q \left(\frac{1}{s^2\va^2}\left( \vert \ppp_tz\vert^2 
+ \sumij \vert \ppp_i\ppp_jz\vert^2\right)
+ \vert \nabla z\vert^2 + s^2\va^2\vert z\vert^2\right)
\weight dxdt\\
\le& C\int_Q \vert z\vert^2 \weight dxdt 
+ C\int_Q \left( \sum_{\vert \gamma\vert\le 2} 
\vert \ppp_x^{\gamma}u\vert^2 + \vert v\vert^2\right) \weight dxdt\\
+ & C\int_Q \vert \ppp_tF\vert^2 \weight dxdt
+ \www{C}(s)(\Vert \ppp_tg\Vert^2_* + \Vert y\Vert_{H^1(\Gamma \times (0,T))}^2)          \\
+ & C\int_{(\ppp\OOO\setminus \Gamma)\times (0,T)}
\left( \frac{\vert \ppp_tg_1\vert^2}{s^2\va^2} + \frac{\vert g_1\vert^2}{\root\of{s\varphi}}
\right) e^{2s\alpha} dSdt
+ C\Vert g_1e^{s\alpha}\Vert^2_{L^2(0,T;H^{\hhalf}(\ppp\OOO))}
+ \www{C}(s)\Vert y\Vert_{H^1(\Gamma \times (0,T))}^2 \\
+& C\int_Q \frac{1}{s\va} \sum_{\vert \gamma\vert\le 2} 
(\vert \ppp_x^{\gamma}u\vert^2 + \vert \ppp_x^{\gamma} v\vert^2) \weight dxdt\\
+& C\int_Q \frac{1}{s\va} \vert \ppp_tG\vert^2 \weight dxdt 
+ \www{C}(s) (\Vert \ppp_th\Vert^2_* + \Vert z\Vert_{H^1(\Gamma \times (0,T))}^2)\\
+ & C\int_{(\ppp\OOO\setminus \Gamma)\times (0,T)}
\left( \frac{\vert \ppp_t(\va^{-\hhalf}h_1)\vert^2}{s^3\va^2} 
+ \frac{s^{-1}\vert \va^{-\hhalf}h_1\vert^2}{\root\of{s\varphi}} \right) e^{2s\alpha} dSdt\\
+ & Cs^{-1}\Vert \va^{-\hhalf}h_1 e^{s\alpha}\Vert^2
_{L^2(0,T;H^{\hhalf}(\ppp\OOO))} 
+ \www{C}(s)\Vert z\Vert_{H^1(\Gamma \times (0,T))}^2 .
\end{align*}
Absorbing the first term on the right-hand side into the left-hand side and 
adding (\ref{(3.5)}) and noting
$$
\frac{1}{s\va} \sum_{\vert \gamma\vert \le 2}\vert \ppp_x^{\gamma}u\vert^2
\le C\sum_{\vert \gamma\vert \le 2}\vert \ppp_x^{\gamma}u\vert^2
$$
in $Q$, we obtain
\begin{align*}
& \int_Q \left(
\frac{1}{s\va}\left( \vert \ppp_ty\vert^2 
+ \sumij \vert \ppp_i\ppp_jy\vert^2\right)
+ s\va\vert \nabla y\vert^2 + s^3\va^3\vert y\vert^2\right) \weight dxdt\\
+ & \int_Q \left(\frac{1}{s^2\va^2}\left( \vert \ppp_tz\vert^2 
+ \sumij \vert \ppp_i\ppp_jz\vert^2\right)
+ \vert \nabla z\vert^2 + s^2\va^2\vert y\vert^2\right)\weight dxdt\\
\le& C\int_Q \left( \vert \ppp_tF\vert^2 + \frac{1}{s\va}\vert \ppp_tG\vert^2
\right) \weight dxdt \\
+ & C\int_Q \left( \sum_{\vert \gamma\vert\le 2} 
\vert \ppp_x^{\gamma}u\vert^2 + \frac{1}{s\va} \sum_{\vert \gamma\vert \le 2}
\vert \ppp_x^{\gamma}v\vert^2 + \vert v\vert^2\right) \weight dxdt\\
+& \www{C}(s)(\Vert \ppp_tg\Vert^2_* + \Vert \ppp_th\Vert^2_*
+ \Vert y\Vert_{H^1(\Gamma \times (0,T))}^2 +\Vert z\Vert_{H^1(\Gamma \times (0,T))}^2 ) \\
+& C\biggl\{ \int_{(\ppp\OOO\setminus \Gamma)\times (0,T)}
\left( \frac{\vert \ppp_tg_1\vert^2}{s^2\va^2} + \frac{\vert g_1\vert^2}{\root\of{s\varphi}}\right)
e^{2s\alpha} dSdt
+ C\Vert g_1e^{s\alpha}\Vert^2_{L^2(0,T;H^{\hhalf}(\ppp\OOO))}\\
+& C\int_{(\ppp\OOO\setminus \Gamma)\times (0,T)}
\left( \frac{\vert \ppp_t(\va^{-\hhalf}h_1)\vert^2}{s^3\va^2} 
+ \frac{s^{-1}\vert \va^{-\hhalf}h_1\vert^2}{\root\of{s\varphi}}\right) e^{2s\alpha} dSdt
+ Cs^{-1} \Vert \va^{-\hhalf}h_1 e^{s\alpha}\Vert^2_{L^2(0,T;
H^{\hhalf}(\ppp\OOO))}\biggr\}.
\end{align*}
Here again we absorb the term $C\int_Q \vert z\vert^2
\weight dxdt$ on the right-hand side, which results from (\ref{(3.5)}),  
into the left-hand side

Applying Theorem \ref{t3} to the second term on the right-hand side, we reach
\begin{eqnarray}\label{(3.7)}
 \int_Q \left(
\frac{1}{s\va}\left( \vert \ppp_ty\vert^2 
+ \sumij \vert \ppp_i\ppp_jy\vert^2\right)
+ s\va\vert \nabla y\vert^2 + s^3\va^3\vert y\vert^2\right) \weight dxdt
                                                       \nonumber\\
+  \int_Q \left(\frac{1}{s^2\va^2}\left( \vert \ppp_tz\vert^2 
+ \sumij \vert \ppp_i\ppp_jz\vert^2\right)
+ \vert \nabla z\vert^2 + s^2\va^2\vert y\vert^2\right)\weight dxdt\nonumber\\
\le C\int_Q \left( s\va\vert F\vert^2 + \vert \ppp_tF\vert^2 + \vert G\vert^2
+ \frac{1}{s\va}\vert \ppp_tG\vert^2 \right) \weight dxdt \nonumber\\
+ \www{C}(s)(\Vert \ppp_tg\Vert^2_* + \Vert \ppp_th\Vert^2_*
 + \Vert g\Vert^2_* + \Vert h\Vert^2_* 
+ \Vert u\Vert_{H^1(\Gamma \times (0,T))}^2 + \Vert v\Vert_{H^1(\Gamma \times (0,T))}^2 \nonumber\\+ \Vert y\Vert_{H^1(\Gamma \times (0,T))}^2 + \Vert z\Vert_{H^1(\Gamma \times (0,T))}^2)\nonumber\\
+ \biggl\{ C\int_{(\ppp\OOO\setminus \Gamma)\times (0,T)}
\left( \frac{\vert \ppp_tg_1\vert^2}{s^2\va^2} + \frac{\vert g_1\vert^2}{\root\of{s\varphi}}\right)
e^{2s\alpha} dSdt
+ C\Vert g_1e^{s\alpha}\Vert^2_{L^2(0,T;H^{\hhalf}(\ppp\OOO))}\biggr\}
                                                    \nonumber\\
+  C\int_{(\ppp\OOO\setminus \Gamma)\times (0,T)}
\left( \frac{\vert \ppp_t(\va^{-\hhalf}h_1)\vert^2}{s^3\va^2} 
+ \frac{s^{-1}\vert \va^{-\hhalf}h_1\vert^2}{\root\of{s\varphi}} \right) e^{2s\alpha} dSdt
\nonumber\nonumber\\
+ Cs^{-1} \Vert \va^{-\hhalf}h_1e^{s\alpha}\Vert^2_{L^2(0,T;
H^{\hhalf}(\ppp\OOO))}\biggr\}. 
\end{eqnarray}
Here we set 
\begin{eqnarray}\label{(3.8)}
I:= \int_{(\ppp\OOO\setminus \Gamma)\times (0,T)}
\left( \frac{\vert \ppp_tg_1\vert^2}{s^2\va^2} + \frac{\vert g_1\vert^2}{\root\of{s\varphi}}\right)
e^{2s\alpha} dSdt\nonumber\\
+  \int_{(\ppp\OOO\setminus \Gamma)\times (0,T)}
\left( \frac{\vert \ppp_t(\va^{-\hhalf}h_1)\vert^2}{s^3\va^2} 
+ \frac{s^{-1}\vert \va^{-\hhalf}h_1\vert^2}{\root\of{s\varphi}} \right) e^{2s\alpha} dSdt\nonumber\\
+  \Vert g_1e^{s\alpha}\Vert^2_{L^2(0,T;H^{\hhalf}(\ppp\OOO))}
+ s^{-1} \Vert \va^{-\hhalf}h_1e^{s\alpha}\Vert^2_{L^2(0,T;
H^{\hhalf}(\ppp\OOO))}
\nonumber\\
=: I_1 + I_2 + I_3 + I_4.
                                         \end{eqnarray}
Now we estimate $I_1, I_2, I_3, I_4$ separately.
\\
{\bf Estimation of $I_1$.}
\\
We can represent 
$$
g_1 = g_{11}(x,t)\cdot \nabla u + g_{10}(x,t)u \quad 
\mbox{on $\ppp\OOO\times (0,T)$},
$$
where $g_{11}, g_{10}$ can be extended to functions in 
$C^1(\ooo{Q})$ by (\ref{lopukh2}).

Therefore, (\ref{(3.1)}) implies
$$
\vert I_1\vert = \left\vert \int_{\ppp\OOO\times (0,T)} 
\left( \frac{1}{s^2\va^2}\vert \ppp_tg_1\vert^2 + \frac{\vert g_1\vert^2}{\root\of{s\varphi}}\right)
 dxdt\right\vert.
$$
Since 
$$
\ppp_tg_1 = g_{11}\cdot \nabla y + g_{10}y + (\ppp_tg_{11})\cdot \nabla u
+ (\ppp_tg_{10})u \quad \mbox{on $\ppp\OOO\times (0,T)$},
$$
we see 
$$
\vert g_1\vert^2 \le C(\vert \nabla u\vert^2 + \vert u\vert^2) 
$$
and
$$
\vert \ppp_tg_1\vert^2 \le C(\vert \nabla y\vert^2 + \vert y\vert^2
+ \vert \nabla u\vert^2 + \vert u\vert^2)
\quad \mbox{on $\ppp\OOO\times (0,T)$}. 
$$
Hence,
$$
\vert I_1\vert \le C\int_{\ppp\OOO\times (0,T)}
\left( \frac{1}{s^2\va^2}(\vert \nabla y\vert^2 + \vert y\vert^2)
+ \frac{1}{\root\of{s\varphi}}(\vert \nabla u\vert^2 + \vert u\vert^2) \right) dSdt.
$$
Consequently (\ref{(3.1)}) and (\ref{(3.2)}) imply
\begin{eqnarray}\label{(3.9)}
\vert I_1\vert 
\le C\int_Q \biggl( \frac{1}{s^2\va^2}\sum_{\vert\gamma\vert\le 2} 
\vert \ppp_x^{\gamma}y\vert^2 
+ \sum_{\vert\gamma\vert\le 2} \vert \ppp_x^{\gamma}u\vert^2
\nonumber\\
+ \vert \nabla y\vert^2 + \vert y\vert^2 + s^2\va^2(\vert \nabla u\vert^2 
+ \vert u\vert^2)\biggr) \weight dxdt.      
                                                \end{eqnarray}
\\
{\bf Estimation of $I_2$}
\\
In view of (\ref{(2.10)}), we have
$$
\vert \ppp_t(\va^{-\hhalf}h_1)\vert = \vert \va^{-\hhalf}\ppp_th_1
 - \hhalf\va^{-\frac{3}{2}}(\ppp_t\va)h_1\vert 
\le C(\va^{-\hhalf}\vert \ppp_th_1\vert + \va^{\hhalf}\vert h_1\vert),
$$
so that 
$$
\vert I_2\vert \le C\int_{\ppp\OOO\times (0,T)}
\left( \frac{1}{s^3\va^3}\vert \ppp_th_1\vert^2 
+ \frac{1}{(s\va)^\frac 32}\vert h_1\vert^2 \right) \weight dSdt.
$$
Since we can represent $h_1 = h_{11}\cdot \nabla v + h_{10}v$, 
it follows that 
$$
\vert h_1\vert^2 \le C(\vert \nabla v\vert^2
+ \vert v\vert^2)
$$
and
$$
\vert \ppp_th_1\vert^2
\le C(\vert \nabla z\vert^2 + \vert z\vert^2 + \vert \nabla v\vert^2
+ \vert v\vert^2) \quad \mbox{on $\ppp\OOO\times (0,T)$}.
$$
Hence, (\ref{(3.1)}) and (\ref{(3.2)}) imply
\begin{eqnarray}\label{(3.10)}
\vert I_2\vert \le C\int_{\ppp\OOO\times (0,T)}
\left( \frac{1}{s^3\va^3}(\vert \nabla z\vert^2 + \vert z\vert^2) 
+ \frac{1}{s\va}(\vert \nabla v\vert^2 + \vert v\vert^2)\right) 
\weight dSdt\nonumber \\
\le C\int_Q \biggl( \frac{1}{s^3\va^3}
\sum_{\vert\gamma\vert\le 2} \vert \ppp_x^{\gamma}z\vert^2
+ \frac{1}{s\va} \sum_{\vert \gamma\vert \le 2} \vert \ppp_x^{\gamma}v\vert^2
\nonumber\\
+ \frac{1}{s\va}(\vert \nabla z\vert^2 + \vert z\vert^2)
+ s\va(\vert \nabla v\vert^2 + \vert v\vert^2)\biggr) \weight dxdt.
\end{eqnarray}
\\
{\bf Estimation of $I_3$}
\\
By noting that $g_1= -\sumij (\ppp_ta_{ij})(\ppp_iu)\nu_j
+ (\ppp_tp)u$ on $\ppp\OOO \times (0,T)$, the trace theorem yields
\begin{eqnarray}\label{(3.11)}
\vert I_3\vert = \Vert g_1e^{s\alpha}\Vert^2_{L^2(0,T;H^{\hhalf}(\ppp\OOO))}
\le C\Vert g_1e^{s\alpha}\Vert^2_{L^2(0,T;H^1(\OOO))}\nonumber\\
= C\int_Q (\vert \nabla(g_1e^{s\alpha})\vert^2 + \vert g_1e^{s\alpha}\vert^2)
dxdt 
\le C\int_Q (\vert \nabla g_1\vert^2 + s^2\va^2\vert g_1\vert^2)\weight dxdt
\nonumber\\
\le C\int_Q \left( \sumij \vert \ppp_i\ppp_j u\vert^2
+ s^2\va^2(\vert \nabla u\vert^2 + \vert u\vert^2) \right)
\weight dxdt.
                                       \end{eqnarray}
\\
{\bf Estimation of $I_4$}
\\
The trace theorem yields
\begin{align*}
& \vert I_4\vert 
= s^{-1}\Vert \va^{-\hhalf}h_1e^{s\alpha}\Vert^2
_{L^2(0,T;H^{\hhalf}(\ppp\OOO))}
\le Cs^{-1}\Vert \va^{-\hhalf}h_1e^{s\alpha}\Vert^2_{L^2(0,T;H^1(\OOO))}\\
\le& Cs^{-1}(\Vert \va^{-\hhalf}h_1e^{s\alpha}\Vert^2_{L^2(Q)}
+ \Vert \nabla(\va^{-\hhalf}h_1e^{s\alpha})\Vert^2_{L^2(Q)})\\
\le& Cs^{-1}\biggl(\Vert \va^{-\hhalf}h_1e^{s\alpha}\Vert^2_{L^2(Q)}
+ \left\Vert \hhalf\va^{-\frac{3}{2}} (\nabla \va)h_1e^{s\alpha}
\right\Vert^2_{L^2(Q)} \\
+ & \Vert s\va^{-\hhalf}h_1(\nabla\va)e^{s\alpha}\Vert^2_{L^2(Q)}
+ \Vert \va^{-\hhalf}(\nabla h_1)e^{s\alpha}\Vert^2_{L^2(Q)}\biggr)\\
\le& C\int_Q \left(s\va\vert h_1\vert^2 + \frac{\vert \nabla h_1\vert^2}
{s\va}\right) \weight dxdt.
\end{align*}
Since $\vert h_1\vert \le C(\vert \nabla v\vert + \vert v\vert)$ and
$\vert \nabla h_1\vert 
\le C\sum_{\vert\gamma\vert\le 2} \vert \ppp_x^{\gamma}v\vert^2$
in $Q$, we have
\begin{equation}\label{(3.12)}
\vert I_4\vert
\le C\int_Q \left( s\va(\vert \nabla v\vert^2 + \vert v\vert^2)
+ \frac{1}{s\va}\sum_{\vert\gamma\vert\le 2} \vert \ppp_x^{\gamma}v\vert^2
\right) \weight dxdt.                        \end{equation}
We rewrite inequality(\ref{(2.17)}) as
\begin{eqnarray}\label{(3.13)}
 \int_Q \biggl\{ \vert \ppp_tu\vert^2 + \sumij \vert \ppp_i\ppp_ju\vert^2
+ s^2\va^2\vert \nabla u\vert^2 + s^4\va^4\vert u\vert^2 
\nonumber\\
+ \frac{1}{s\va}\left(
\vert \ppp_tv\vert^2 + \sumij \vert \ppp_i\ppp_jv\vert^2\right)
+ s\va\vert \nabla v\vert^2 + s^3\va^3\vert v\vert^2 \biggr\}
\weight dxdt \le J,                
\end{eqnarray}
where we set 
$$
 J := C\int_Q (s\va\vert F\vert^2 + \vert G\vert^2) \weight dxdt
+ \www{C}(s)(\Vert g\Vert_*^2 + \Vert h\Vert^2_* +\Vert u\Vert_{H^1(\Gamma \times (0,T))}^2+\Vert v\Vert_{H^1(\Gamma \times (0,T))}^2 ).
$$
Applying (\ref{(3.13)}) to (\ref{(3.9)}) - (\ref{(3.12)}), we have
\begin{equation}\label{(3.14)}
\left\{ \begin{array}{rl}
& \vert I_1\vert \le J 
+ C\int_Q \left( \frac{1}{s^2\va^2}\sum_{\vert \gamma\vert\le 2}
\vert \ppp_x^{\gamma}y\vert^2 + \vert \nabla y\vert^2 + \vert y\vert^2
\right) \weight dxdt, \cr\\
& \vert I_2\vert \le J 
+ C\int_Q \left( \frac{1}{s^3\va^3}\sum_{\vert \gamma\vert\le 2}
\vert \ppp_x^{\gamma}z\vert^2 
+ \frac{1}{s\va}(\vert \nabla z\vert^2 + \vert z\vert^2)
\right) \weight dxdt,\cr \\
& \vert I_3\vert\le J\quad \mbox{and} \quad \vert I_4\vert \le J.
\end{array}\right.
\end{equation}

Applying (\ref{(3.14)}) to (\ref{(3.7)}), we reach
\begin{align*}
& \int_Q \left(
\frac{1}{s\va}\left( \vert \ppp_ty\vert^2 
+ \sumij \vert \ppp_i\ppp_jy\vert^2\right)
+ s\va\vert \nabla y\vert^2 + s^3\va^3\vert y\vert^2\right) \weight dxdt\\
+ & \int_Q \left(\frac{1}{s^2\va^2}\left( \vert \ppp_tz\vert^2 
+ \sumij \vert \ppp_i\ppp_jz\vert^2\right)
+ \vert \nabla z\vert^2 + s^2\va^2\vert y\vert^2\right)\weight dxdt\\
\le& C\int_Q \left( s\va\vert F\vert^2 + \vert \ppp_tF\vert^2 + \vert G\vert^2
+ \frac{1}{s\va}\vert \ppp_tG\vert^2 \right) \weight dxdt \\
+ & \www{C}(s)(\Vert \ppp_tg\Vert^2_* + \Vert g\Vert^2_* 
+ \Vert \ppp_th\Vert^2_* + \Vert h\Vert^2_* \nonumber\\
+& \Vert u\Vert_{H^1(\Gamma \times (0,T))}^2 + \Vert v\Vert_{H^1(\Gamma \times (0,T))}^2 + \Vert y\Vert_{H^1(\Gamma \times (0,T))}^2 + \Vert z\Vert_{H^1(\Gamma \times (0,T))}^2) + J\\
+& \int_Q \left( \frac{1}{s^2\va^2}\sum_{\vert \gamma\vert\le 2}
\vert \ppp_x^{\gamma}y\vert^2 + \vert \nabla y\vert^2 + \vert y\vert^2
+ \frac{1}{s^3\va^3}\sum_{\vert \gamma\vert\le 2}
\vert \ppp_x^{\gamma}z\vert^2 + \frac{1}{s\va}(\vert \nabla z\vert^2 
+ \vert z\vert^2)\right) \weight dxdt
\end{align*}
for all large $s>0$.
Choosing $s>0$ large, we can absorb the final term on the right-hand side into
the left-hand side, we complete the proof of Proposition \ref{p1}.
$\blacksquare$
\\
{\bf Third Step: Completion of the proof.}
\\
Since $F(x,t)= q_1(x,t)f_1(x)$ and $G(x,t) = q_2(x,t)f_2(x,t)$,
by Proposition \ref{p1}, we have
\begin{eqnarray}\label{(3.15)}
\int_Q \left(
\frac{1}{s\va}\vert \ppp_t^2u\vert^2 + s^3\va^3\vert \ppp_tu\vert^2
+ \frac{1}{s^2\va^2} \vert \ppp_t^2v\vert^2 
+ s^2\va^2 \vert \ppp_tv\vert^2\right) \weight dxdt
\nonumber\\
\le C\int_Q (s\va \vert f_1\vert^2 + \vert f_2\vert^2) \weight dxdt
+ \www{C}(s)J_0                         
\end{eqnarray}
for all large $s>0$.  Here we set 
$$
J_0:= \sum_{k=0}^1 (\Vert \ppp_t^kg\Vert^2_* 
+ \Vert \ppp_t^kh\Vert^2_* +\Vert \ppp_t^ku\Vert_{H^1(\Gamma \times (0,T))}^2+\Vert \ppp_t^kv\Vert_{H^1(\Gamma \times (0,T))}^2 ).
$$

We note that (\ref{(2.10)}) and 
$\lim_{t\to 0} \alpha(x,t) = -\infty$ for $x\in \OOO$.
Then, 
\begin{align*}
& \int_{\OOO} s\va(x,t_0) \vert \ppp_tu(x,t_0)\vert^2
e^{2s\alpha(x,t_0)} dx 
= \int^{t_0}_0 \ppp_t\left( \int_{\OOO} s\va\vert \ppp_tu(x,t)\vert^2
e^{2s\alpha(x,t)} dx \right) dt \\
=& \int^{t_0}_0 \int_{\OOO} (2s\va (\ppp_tu)(\ppp_t^2u)
+ s(\ppp_t\va)\vert \ppp_tu\vert^2 
+ 2s^2\va (\ppp_t\alpha)\vert \ppp_tu\vert^2) \weight dxdt\\
\le& C\int_Q (s\va\vert \ppp_tu\vert \vert \ppp_t^2u\vert
+ s\va^2 \vert \ppp_tu\vert^2 + s^2\va^3\vert \ppp_tu\vert^2) \weight dxdt \\
= & C\int_Q \{ (s^{-\hhalf}\va^{-\hhalf}\vert \ppp_t^2u\vert)
(s^{\frac{3}{2}}\va^{\frac{3}{2}} \vert \ppp_tu\vert) 
+ (s\va^2 + s^2\va^3) \vert \ppp_tu\vert^2\} \weight dxdt \\
\le & C\int_Q \left(\frac{1}{s\va}\vert \ppp_t^2u\vert^2
+ (s^3\va^3 + s\va^2 + s^2\va^3)\vert \ppp_tu\vert^2\right) \weight dxdt\\
\le & C\int_Q \left(\frac{1}{s\va}\vert \ppp_t^2u\vert^2 
+ s^3\va^3\vert \ppp_tu\vert^2\right) \weight dxdt.
\end{align*}
For the last inequality, we used $s\va^2 + s^2\va^3 \le Cs^3\va^3$ for
large $s>0$.
Hence, by (\ref{(3.15)}), we have
\begin{equation}\label{(3.16)}
\int_{\OOO} s\va(x,t_0) \vert \ppp_tu(x,t_0)\vert^2
e^{2s\alpha(x,t_0)} dx 
\le C\int_Q \left(s\va \vert f_1\vert^2 + \vert f_2\vert^2\right) 
\weight dxdt + \www{C}(s)J_0
\end{equation}
for all large $s>0$.

Next, by (\ref{(3.15)}), we can similarly estimate
\begin{eqnarray}\label{(3.17)}
 \int_{\OOO} \vert \ppp_tv(x,t_0)\vert^2
e^{2s\alpha(x,t_0)} dx 
= \int^{t_0}_0 \ppp_t\left( \int_{\OOO} \vert \ppp_tv(x,t)\vert^2
e^{2s\alpha(x,t)} dx \right) dt \nonumber\\
= \int^{t_0}_0 \int_{\OOO} (2 (\ppp_tv)(\ppp_t^2v)
+ 2s(\ppp_t\alpha)\vert \ppp_tv\vert^2) \weight dxdt\nonumber\\
\le C\int_Q (\vert \ppp_tv\vert \vert \ppp_t^2v\vert
+ s\va^2 \vert \ppp_tv\vert^2) \weight dxdt \nonumber\\
= C\int_Q \left\{
\left( \frac{1}{s\va} \vert \ppp_t^2v\vert\right) (s\va \vert \ppp_tv\vert) 
+ s\va^2 \vert \ppp_tv\vert^2 \right\} \weight dxdt  \nonumber           \\
\le  C\int_Q \left( \frac{1}{s^2\va^2} \vert \ppp_t^2v\vert^2
+ s^2\va^2\vert \ppp_tv\vert^2 + s\va^2\vert \ppp_tv\vert^2
\right) \weight dxdt
\nonumber\\
\le  C\int_Q (s\va \vert f_1\vert^2 + \vert f_2\vert^2) \weight dxdt 
+ \www{C}(s)J_0
\end{eqnarray}
for all large $s>0$.

On the other hand, by assumption (\ref{(1.5)}) we have
$$
\left\{ \begin{array}{rl}
& f_1(x) = \frac{1}{q_1(x,t_0)}\ppp_tu(x,t_0)
+ \frac{1}{q_1(x,t_0)}(A(t_0)u_0 - c_0(x,t_0)v_0),\cr\\
& f_2(x) = \frac{1}{q_2(x,t_0)}\ppp_tv(x,t_0)
+ \frac{1}{q_2(x,t_0)}(-B(t_0)v_0 - A_0(t_0)u_0), \quad x\in \OOO
\end{array}\right.
$$
and
$$
\left\{ \begin{array}{rl}
& s\va(x,t_0)\vert f_1(x)\vert^2
\le Cs\va(x,t_0)\vert \ppp_tu(x,t_0)\vert^2
+ Cs\va(x,t_0) \left( \sum_{\vert \gamma\vert\le 2}
\vert \ppp_x^{\gamma}u_0(x)\vert^2 + \vert v_0(x)\vert^2\right),\\
& \vert f_2(x)\vert^2
\le C\vert \ppp_tv(x,t_0)\vert^2
+ C\sum_{\vert \gamma\vert\le 2}(\vert \ppp_x^{\gamma}u_0(x)\vert^2 
+ \vert \ppp_x^{\gamma} v_0(x)\vert^2), \quad x\in \OOO.
\end{array}\right.
$$
Hence (\ref{(3.16)}) and (\ref{(3.17)}) yield
\begin{eqnarray}\label{(3.18)}
 \int_{\OOO} (s\va(x_0,t)\vert f_1(x)\vert^2 + \vert f_2(x)\vert^2)
e^{2s\alpha(x,t_0)} dx\nonumber\\
\le C\int_{\OOO} (s\va(x,t_0)\vert \ppp_tu(x,t_0)\vert^2
+ \vert \ppp_tv(x,t_0)\vert^2) \weight dxdt
+ \www{C}(s)(\Vert u_0\Vert^2_{H^2(\OOO)} + \Vert v_0\Vert^2_{H^2(\OOO)})
\nonumber\\
\le C\int_Q \left( s\va\vert f_1\vert^2 + \vert f_2\vert^2\right)
\weight dxdt + \www{C}(s)J_1  
\end{eqnarray}              
for all large $s>0$.
Here and henceforth we set
$$
J_1:= J_0 + \Vert u_0\Vert^2_{H^2(\OOO)}
  + \Vert v_0\Vert^2_{H^2(\OOO)}.
$$
On the other hand, we see
\begin{align*}
& \va(x,t)\vert f_1(x)\vert^2 e^{2s\alpha(x,t)}
= \va(x,t_0)\vert f_1(x)\vert^2 e^{2s\alpha(x,t_0)}
\times \frac{\va(x,t)}{\va(x,t_0)}e^{2s(\alpha(x,t)-\alpha(x,t_0))}\\
=& \va(x,t_0)\vert f_1(x)\vert^2 e^{2s\alpha(x,t_0)}
\times \frac{\mu(t_0)}{\mu(t)}
e^{-2s \xi(x)\left( \frac{1}{\mu(t)} - \frac{1}{\mu(t_0)}\right)}.
\end{align*}
Here and henceforth we set 
$$
\xi(x):= e^{2\la\Vert {\eta}\Vert_{C(\ooo{\OOO})}}
- e^{\la{\eta}(x)}> 0, \quad
C_1:= e^{2\la\Vert {\eta}\Vert_{C(\ooo{\OOO})}}
- e^{\la\Vert {\eta}\Vert_{C(\ooo{\OOO})}} 
= \min_{x\in \ooo{\OOO}}\xi(x) > 0.
$$
Therefore,
\begin{align*}
& \int_Q \va(x,t)\vert f_1(x)\vert^2 e^{2s\alpha(x,t)} dxdt
\le C\int_Q \va(x,t_0)\vert f_1(x)\vert^2 e^{2s\alpha(x,t_0)}
\frac{1}{\mu(t)}
e^{-2s\xi(x)\left( \frac{1}{\mu(t)} - \frac{1}{\mu(t_0)}\right)} dxdt\\
\le& C\int_{\OOO} \va(x,t_0)\vert f_1(x)\vert^2 e^{2s\alpha(x,t_0)}
\left( \int^T_0 \frac{1}{\mu(t)}
e^{-2sC_1\left( \frac{1}{\mu(t)} - \frac{1}{\mu(t_0)}\right)}dt
\right) dx.
\end{align*}
We will estimate 
$\int^T_0 \frac{1}{\mu(t)}
e^{-2sC_1\left( \frac{1}{\mu(t)} - \frac{1}{\mu(t_0)}\right)}dt$.
Indeed 
$$
\lim_{s\to \infty} 
\frac{1}{\mu(t)}
e^{-2sC_1\left( \frac{1}{\mu(t)} - \frac{1}{\mu(t_0)}\right)} = 0
$$
for each fixed $t \in (0,T) \setminus \{t_0\}$.  Next, since 
$$
\frac{1}{\mu(t)}e^{-2sC_1\left( \frac{1}{\mu(t)} - \frac{1}{\mu(t_0)}\right)}
\le \frac{1}{\mu(t)}
e^{-2C_1\left( \frac{1}{\mu(t)} - \frac{1}{\mu(t_0)}\right)}
$$
for $s \ge 1$ and by $\mu(t) \le \mu(t_0)$, we have
\begin{align*}
& \sup_{s\ge 1}\sup_{0< t< T} \frac{1}{\mu(t)}
\exp\left( -2sC_1\left( \frac{1}{\mu(t)} - \frac{1}{\mu(t_0)} \right)\right)
\le \sup_{0<t<T} \left( \frac{1}{\mu(t)}
e^{-\frac{2C_1}{\mu(t)}}\right) e^{\frac{2C_1}{\mu(t_0)}}\\ 
\le &\sup_{\zeta>0} (\zeta e^{-2C_1\zeta}) e^{\frac{2C_1}{\mu(t_0)}} < \infty.
\end{align*}
Consequently, the Lebesgue convergence theorem yields
$$
\int^T_0 \frac{1}{\mu(t)}
e^{-2sC_1\left( \frac{1}{\mu(t)} - \frac{1}{\mu(t_0)}\right)} dt
= o(1) \quad \mbox{as $s \to \infty$}.
$$
Hence,
$$
\int_Q \va(x,t)\vert f_1(x)\vert^2 e^{2s\alpha(x,t)} dxdt
= o(1)\int_{\OOO} \va(x,t_0)\vert f_1(x)\vert^2 e^{2s\alpha(x,t_0)} dx
$$
as $s \to \infty$.  

Similarly we can verify 
$$
\int_Q \vert f_2(x)\vert^2 e^{2s\alpha(x,t)} dxdt
= o(1)\int_{\OOO} \vert f_2(x)\vert^2 e^{2s\alpha(x,t_0)} dx
$$
as $s \to \infty$.  
Therefore, (\ref{(3.18)}) yields
\begin{align*}
& \int_{\OOO} (s\va(x,t_0)\vert f_1(x)\vert^2 + \vert f_2(x)\vert^2)
e^{2s\alpha(x,t_0)} dx\\
= &o(1)\int_{\OOO} (s\va(x,t_0)\vert f_1(x)\vert^2 + \vert f_2(x)\vert^2)
e^{2s\alpha(x,t_0)} dx
+ \www{C}(s)J_1
\end{align*}
for all large $s>0$.  Choosing $s>0$ large, we can absorb the first term on 
the right-hand side into the left-hand side to obtain
$$
\int_{\OOO} (s\va(x,t_0)\vert f_1(x)\vert^2 + \vert f_2(x)\vert^2)
e^{2s\alpha(x,t_0)} dx
\le \www{C}(s)J_1.
$$
With such fixed $s>0$, we have 
$\min_{x\in \ooo{\OOO}} s\va(x,t_0)e^{2s\alpha(x,t_0)} > 0$, and so
we complete the proof of Theorem \ref{t2}.
$\blacksquare$
\section{State determination for the nonlinear mean field game system}
\label{Z4}
For the non-linearized mean field game system (\ref{(1.1)}), we discuss the 
state determination problem.  In this article, we do not consider inverse 
source problems and inverse coefficient problems of determining spatially 
varying factors of the coefficients such as $\kappa(x,t)$, 
and we postpone them to a future work.  

In (\ref{(1.1)}), we assume that 
\begin{eqnarray}\label{(4.1)}
\left\{ \begin{array}{rl}
 a, \, \nabla a \in C^1(\ooo{Q}), \quad a > 0 \quad 
\mbox{on $\ooo{Q}$}, \\
 c_0\in L^{\infty}(Q), \quad \kappa \in L^{\infty}(0,T;W^{1,\infty}(\OOO)).
\end{array}\right.
\end{eqnarray}
Then we can prove
\\
\begin{thm}\label{t4} Let (\ref{lopukh}) and (\ref{(4.1)}) holds true.
{\it
For $k=1,2$, let $(u_k,v_k) \in H^{2,1}(Q)$ satisfy  the system of equations
\begin{equation}\label{(4.2)}
\left\{ \begin{array}{rl}
& \ppp_tu_k(x,t) + a(x,t)\Delta u_k(x,t) 
- \hhalf\kappa(x,t)\vert \nabla u_k(x,t)\vert^2 + c_0(x,t)v_k(x,t) = F_k(x,t), 
                                                   \\
& \ppp_tv_k(x,t) - \Delta (a(x,t) v_k(x,t)) 
- \ddd(\kappa(x,t)v_k(x,t)\nabla u_k(x,t)) = G_k(x,t)
\quad \mbox{in $Q$,}
\end{array}\right.
\end{equation}
and
$$
\nabla u_k\cdot \nu = \nabla (av_k) \cdot \nu = 0 
\quad \mbox{on $\ppp\OOO\times (0,T)$}.
$$
We further assume
\begin{equation}\label{(4.2)}
\Vert u_k\Vert_{L^{\infty}(0,T;W^{2,\infty}(\OOO))}
+ \Vert v_k\Vert_{L^{\infty}(0,T;W^{1,\infty}(\OOO))} \le M_1,  
                                             \end{equation}
where $M_1>0$ is an arbitrarily chosen constant.
Then, for arbitrarily given $\ep > 0$, we can find a constant $C_{\ep} >0$ 
such that 
\begin{align*}
& \Vert u_1-u_2\Vert_{\HHHH(\OOO\times (\ep,T-\ep))}
+ \Vert v_1-v_2\Vert_{\HHHH(\OOO\times (\ep,T-\ep))}\\
\le & C_{\ep}(\Vert F_1-F_2\Vert_{L^2(Q)} + \Vert G_1-G_2\Vert_{L^2(Q)}
+ \Vert u_1-u_2\Vert_{H^1(\Gamma\times (0,T))}
+ \Vert v_1-v_2\Vert_{H^1(\Gamma\times (0,T))}).
\end{align*}
Here the constant $C_{\ep}>0$ depends on $\ep, M_1 > 0$ and 
the coefficients $a, c_0, \kappa$.
}
\end{thm}
{\bf Proof of Theorem \ref{t4}.}
\\
We take the difference $u:= u_1-u_2$ and $v:= v_1 - v_2$.  Then 
we note that $\nabla (av_k) = 0$ on $\ppp\OOO$ means 
$$
\ppp_{\nu}v_k + \frac{\ppp_{\nu}a}{a}v_k = 0 \quad 
\mbox{on $\ppp\OOO \times (0,T)$, $k=1,2$.}
$$
Hence,
$$
\left\{ \begin{array}{rl}
& \ppp_tu + a(x,t)\Delta u(x,t) 
- \hhalf\kappa(x,t)((\nabla u_1 + \nabla u_2)\cdot\nabla u) + c_0(x,t)v 
= F_1-F_2, \\
& \ppp_tv - a(x,t)\Delta v(x,t) -2\nabla a(x,t)\cdot \nabla v(x,t)
- v(x,t)\Delta a(x,t)\\ 
- & (\kappa\nabla u_1\cdot \nabla v 
+ (\nabla\kappa \cdot \nabla u_1 + \kappa\Delta u_1)v 
+ v_2\kappa \Delta u + \nabla(\kappa v_2)\cdot \nabla u)
= G_1 - G_2, \quad (x,t) \in Q, \\
& \ppp_{\nu} u = \ppp_{\nu}v + \frac{\ppp_{\nu}a}{a} v = 0 \quad 
\mbox{on $\ppp\OOO\times (0,T)$}.
\end{array}\right.
$$

By (\ref{(4.1)}) and (\ref{(4.2)}), we can find a constant $M_2 > 0$ depending 
on 
$M_1$ such that 
\begin{align*}
& \left\Vert \sum_{k=1}^2 \kappa\nabla u_k\right\Vert_{L^{\infty}(Q)}
+ \Vert\nabla \kappa \cdot \nabla u_1 + \kappa\Delta u_1\Vert_{L^{\infty}(Q)}\\
+& \Vert \kappa v_2\Vert_{L^{\infty}(Q)}  
+ \Vert\nabla (\kappa v_2)\Vert_{L^{\infty}(Q)} \le M_2
\end{align*}
and
$$
\frac{\ppp_{\nu}a}{a} \in C^1(\ppp\OOO \times [0,T]).
$$
Setting $p:= 0$ and $q:= -\frac{\ppp_{\nu}a}{a}$ in (\ref{(1.4)}), 
we can satisfy the conditions of Theorem \ref{t1} in view of 
(\ref{(4.1)}).  
Thus the proof of Theorem \ref{t4} follows directly from Theorem \ref{t1}.
$\blacksquare$

\section{Proof of Lemma \ref{1}}\label{Z5}

We set $u(x,t) = w(x,T-t)$ for $(x,t) \in Q$.
Then, since $\alpha(x,t) = \alpha(x,T-t)$ and
$\va(x,t) = \va(x,T-t)$ for $(x,t) \in Q$, the change of variables 
$t \mapsto T-t$ transfer the Carleman estimate for $\ppp_t + A(t)$ to 
$\ppp_t - A(t)$.
Thus it is sufficient to prove Lemma \ref{1} for the parabolic operator 
$\ppp_t - \sumij a_{ij}\ppp_i\ppp_j$, which is forward in time.

Indeed we will prove a sharper estimate than Lemma \ref{1}.
In order to formulate our estimate we introduce the operators 
\begin{equation}\label{(5.9')}
\left\{ \begin{array}{rl}
& L_2(x,t,D,s) w  = - \sumij a_{ij} \ppp_i\ppp_jw - s^2\lambda^2\varphi^2 
a(x,t, \nabla\psi, \nabla\psi) w -  s(\ppp_t\alpha) w,  \\
& L_1(x,t,D,s) w  = \partial_t w + 2s\lambda \varphi
\sumij a_{ij} (\ppp_i\psi)\ppp_jw + 2s\lambda^2 \varphi
a(x,t, \nabla\psi, \nabla\psi) w.
\end{array}\right.
\end{equation}
We have
\begin{lem}\label{6}
{\it 
Let $F\in L^2(Q), g\in L^2(0,T;H^\frac 12(\partial\Omega)), 
\partial_t g\in L^2((\partial\Omega\setminus\Gamma)\times (0,T))$ and 
$p\in C^1(\partial\Omega\times [0,T])$.
There exists a constant $\lambda_0 > 0$ such that for an arbitrary 
$\lambda \ge \la_0$, we can choose a constant $s_0(\lambda)$ 
satisfying : there exists a constant $C>0$ such that 
\begin{align*}
& \int_Q \left( \frac{1}{s\varphi} \left(\vert\partial_t u\vert^2 
+ \sumij \vert \ppp_i\ppp_ju\vert^2\right) 
+ s\lambda^2 \varphi |\nabla u|^2 + s^3 \lambda^4\varphi^3\vert u\vert^2\right)
e^{2s\alpha}  dx\,dt \\
+ & \sum_{k=1}^2\Vert L_k(x,t,D,s) (ue^{{s}\alpha})\Vert^2_{L^2(Q)}  \\
\le &C\biggl( 
\int_Q |F|^2 e^{2s\alpha}  dxdt 
+ \int_{(\ppp\OOO\setminus \Gamma)\times (0,T)}
\left(\frac{\vert \partial_t g\vert^2}{s^2\lambda^2\varphi^2}
+ \frac{1}{\root\of{s\varphi}} \vert g\vert^2\right) e^{2s\alpha} dSdt 
+ \Vert ge^{s\alpha}\Vert^2_{L^2(0,T;H^\frac 12(\partial\Omega))} \\
+& \int_{\Gamma \times (0,T)} 
\left( s\lambda \varphi \vert \nabla u\vert^2  
+ s^3\lambda^3\varphi^3 \vert u\vert^2 
+ \frac{\vert \partial_t u\vert^2}{s\varphi} \right) e^{2s\alpha} dS dt
\biggr).
\end{align*}
for all $s> s_0(\la)$ and all $u\in H^{2,1}(Q)$ satisfying
$$
\left\{ \begin{array}{rl}
& \ppp_tu(x,t) - \sumij a_{ij}(x,t)\ppp_i\ppp_ju(x,t) = F(x,t), \quad
(x,t) \in Q, \\
& \ppp_{\nu_A} u(x,t) - p(x,t)u(x,t) = g(x,t), \quad (x,t) \in \ppp\OOO\times
(0,T).
\end{array}\right.
$$ 
}
\end{lem}

Fixing $\la>0$ sufficiently large, we obtain Lemma \ref{1} 
from Lemma \ref{6}. We recall that the bounds $M, M_0$ are defined by 
(\ref{(2.2)}).

Section 5 is now devoted 
\\
{\bf Proof of Lemma \ref{6}.}
\\
{\bf First Step.}
\\
Let us introduce a quadratic form:
$$
a(x,t,V,W) : =\sumij a_{ij}(x,t)v_iw_j \quad \mbox{for 
$V:= (v_1,..., v_d)$ and $W:= (w_1, ..., w_d)$}.    
$$

We recall that $\Gamma_1$ is a relatively open subboundary of $\ppp\OOO$ and 
$\ooo{\ppp\OOO \setminus \Gamma} \subset \Gamma_1$.
Let $\mathcal U$ be a subdomain of $\Omega$ such that 
$\mathcal{U}\cap \partial\Omega \subset \Gamma_1$.
Without loss of generality, we can assume that
$$
\mbox{supp} \,u \subset \mathcal U\times [0,T].
$$ 
Indeed, let $\mathcal U_1$ be an open set such that 
$\Omega\subset \mathcal U\cup\mathcal U_1$ and $\overline{\mathcal U_1}\cap 
\ooo\Gamma =\emptyset$.  Let $e_1, e_2\in C^\infty_0(\Bbb R^n)$ 
be a partition of unity subject to the covering $\mathcal U_1,\mathcal U_2$.
Similarly to \cite {Hor}, it suffices to prove the Carleman estimate 
Lemma (\ref{6}) for the functions $ue_1$ and $ue_2$.
The proof for the function $e_2u$ is simpler, since it does not require 
consideration of the function $\widetilde w$ (introduced below in (\ref{im})) and follows directly from 
the inequality (\ref{(5.24)}) derived below.

We consider the operator estimate 
$$
 \widehat{L}(x,t,D)u = \partial_t u - \sumij a_{ij}(x,t)\ppp_i\ppp_ju
$$ 
and
\begin{equation}\label{(5.1)}
\widehat{L}(x,t,D) u = \widetilde F \quad \mbox{in $Q$},           
\end{equation}
where 
\begin{equation}\label{(5.2)}
\www{F}(x,t) = F(x,t) +  \sumij (\ppp_ia_{ij})(x,t) \ppp_ju.
                      \end{equation}

We set $\widetilde \psi(x) = -\psi(x)$ in a neighborhood $U$. 
Using the function $\widetilde \psi$, we introduce functions 
$\widetilde\alpha$ and $\widetilde\varphi:$
\begin{equation}\label{(5.3)}
  \widetilde\alpha(x,t) = \frac{e^{\lambda\widetilde \psi(x)}
- e^{2\lambda ||\psi||_{C(\overline{\Omega})}}}{\mu(t)},
\quad \widetilde \varphi(x,t) = \frac{e^{\lambda\widetilde \psi(x)}}{\mu(t)}.
\end{equation}

We denote 
\begin{equation}\label{im} w(x,t)=e^{s\alpha}u(x,t)\quad \mbox{and} 
\quad \widetilde w(x,t) = e^{s\widetilde \alpha}u(x,t).\end{equation}

Then, we have
\begin{equation}\label{(5.4)}
w(\cdot,0)=w(\cdot,T) = \widetilde w(\cdot,0) = \widetilde w(\cdot,T)= 0
\quad \mbox{in}\quad \Omega.
\end{equation}
We define operators $P(x,t,D,s)$ and $\widetilde P(x,t,D,s)$ by 
\begin{equation}\label{(5.5)}
P(x,t,D,s)w = e^{s\alpha} \widehat{L}(x,t,D) e^{-s\alpha} w,\quad
\widetilde{P}(x,t,D,s) w = e^{s\widetilde{\alpha}} \widehat{L}(x,t,D) 
e^{-s\widetilde{\alpha}} w.             \end{equation}

It follows from (\ref{(5.1)}) and (\ref{(5.2)}) that 
\begin{equation}\label{(5.6)}
P(x,t,D,s)w  = e^{s\alpha} \widehat{L}(x,t,D)(e^{-s\alpha} w)
= e^{s\alpha} \widetilde{F} \quad \mbox{in}\quad  Q,
                                           \end{equation}
and
\begin{equation}\label{(5.7)}
\widetilde{P}(x,t,D,s)\widetilde{w} = e^{s\widetilde{\alpha}} \widehat{L}(x,t,D) 
e^{-s\widetilde{\alpha}} \widetilde{w}
= e^{s\widetilde{\alpha}} \widetilde{F} \quad \mbox{in}\quad  Q.
                                               \end{equation}

The operator $P$ can be written explicitly as follows
\begin{eqnarray}\label{(5.8)}
P(x,t,D,s)w = \partial_t w - \sumij a_{ij}\ppp_i\ppp_jw 
+ 2\lambda\varphi \sumij a_{ij} (\ppp_i\psi)\ppp_jw
  + s\lambda^2 \varphi a(x,t, \nabla\psi, \nabla\psi) w
\nonumber\\
 - s^2\lambda^2\varphi^2  a(x,t, \nabla\psi, \nabla\psi) w
 + s\lambda \varphi w \sumij a_{ij} \ppp_i\ppp_j\psi - s(\ppp_t\alpha) w.
                                             \end{eqnarray}
We introduce the operators $\www{L}_k(x,t,D,s)$, $k=1,2$ 
as follows.
\begin{equation}\label{(5.9)}
\left\{ \begin{array}{rl}
& \widetilde L_2(x,t,D,s) w  = - \sumij a_{ij} \ppp_i\ppp_jw 
- s^2\lambda^2\widetilde \varphi^2 a(x,t, \nabla\www\psi, \nabla\www\psi) w 
- s(\ppp_t\widetilde \alpha) w, \\
& \widetilde L_1(x,t,D,s) w  = \partial_t w + 2{{s}}\lambda\widetilde\varphi
\sumij a_{ij} (\ppp_i\widetilde\psi) \ppp_jw + 2\lambda^2 \widetilde \varphi
a(x,t, \nabla\widetilde \psi, \nabla\widetilde \psi) w.
\end{array}\right.
                                            \end{equation}
Then,
\begin{equation}\label{(5.10)}
L_1(x,t,D,s) w + L_2(x,t,D,s)w= H(x,t,\lambda,{s})\quad\mbox{in}\,\,Q,
                                  \end{equation}
where 
\begin{equation}\label{(5.11)}
H(x,t,\lambda,s) := \widetilde ge^{s\alpha}+2s\lambda^2 \varphi
a(x,t, \nabla\psi, \nabla\psi) w - s\lambda \varphi w \sumij a_{ij}
(\ppp_i\ppp_j\psi) w                        
\end{equation}
and
\begin{equation}\label{(5.12)}
\widetilde L_1(x,t,D,s) \widetilde w 
+ \widetilde L_2(x,t,D,s)\widetilde w
= \widetilde{H}(x,t,\lambda,s) \quad\mbox{in}\,\,Q,
                                         \end{equation}
where 
\begin{equation}\label{(5.13)}
\widetilde{H}(x,t,\lambda,s) = \widetilde ge^{s\alpha} + 2s\lambda^2 \widetilde\varphi
a(x,t, \nabla\widetilde\psi, \nabla\widetilde \psi) w
- s\lambda \widetilde\varphi\widetilde w 
\sumij a_{ij}(\ppp_i\ppp_j\widetilde  \psi)\widetilde w.
\end{equation}
\\
Henceforth we set 
$$
\Sigma:= \ppp\OOO \times (0,T), \quad d\Sigma := dS dt.
$$
\\
{\bf Second Step.}
\\ 
We will verify the following equality:
\begin{eqnarray}\label{5.14}
 (L_2 w, L_1 w)_{L^2 (Q)}
= \int_Q L_1(x,t,D,s)w \sumij (\ppp_ja_{ij})(\ppp_iw) dxdt \nonumber\\
+ \int_Q \biggl\{
 -\sumij \frac 12 (\ppp_t a_{ij})(\ppp_iw)(\ppp_jw) 
+ \partial_t\left( \frac{s^2\lambda^2 \varphi^2}{2}
a(x,t, \nabla\psi, \nabla\psi)\right) w^2
  + \frac{s\ppp_t^2\alpha}{2} w^2\nonumber\\
+  s^3 \la^4\varphi^3 a(x,t, \nabla\psi, \nabla\psi)^2 w^2
- 2s^2\lambda^2\va (\ppp_t\alpha) a(x,t, \nabla\psi, \nabla\psi) w^2\nonumber\\
+  s\lambda^2 \varphi a(x,t, \nabla\psi, \nabla\psi)a(x,t,\nabla w,\nabla w)\nonumber\\
+  2s\lambda^2 w \sumij a_{ij} (\ppp_jw)
\ppp_i(\varphi a(x,t, \nabla\psi, \nabla\psi)) 
+ 2s\lambda^2 \varphi a(x,t, \nabla\psi, \nabla w)^2\nonumber\\
+ 2s\lambda \varphi \sumij  a_{ij}\ppp_iw
\left(\sum_{k,\ell=1}^d \ppp_j(a_{k\ell}\ppp_k\psi) \ppp_\ell w \right)
- s\lambda\varphi \sum_{k,\ell=1}^d a_{k\ell} (\ppp_k\psi) 
\left( \sumij (\ppp_\ell a_{ij})(\ppp_i w) \ppp_jw\right)   \nonumber         \\
-  s\lambda\varphi \sum_{k,\ell=1}^d a_{k\ell} (\ppp_k\psi)
 \left( \sumij (\ppp_\ell a_{ij})(\ppp_iw)\ppp_jw\right) 
 - a(x,t, \nabla w, \nabla w) s\lambda\varphi \sum_{k,\ell=1}^d
\ppp_\ell(a_{k\ell}\ppp_k\psi) \biggr\} dxdt \nonumber   \\
+  \int_\Sigma (2s\lambda\varphi a(x,t,\nu,\nabla w)a(x,t,\nabla\psi,\nabla w)
 - s\lambda\varphi a(x,t, \nabla w, \nabla w)a(x,t, \nu,\nabla \psi)) 
                                                      d\Sigma\nonumber\\ 
- \int_\Sigma (2s^3\lambda^3\varphi^3 a(x,t, \nabla\psi, \nabla\psi)
a(x,t, \nabla\psi, \nu)   + 2s^2 \lambda\va (\ppp_t\alpha) 
a(x,t, \nabla\psi, \nu))w^2 d\Sigma 
\nonumber\\
- \int_\Sigma a(x,t, \nu,\nabla w)(\partial_t w 
+ 2s\lambda^2\varphi a(x,t, \nabla\psi, \nabla\psi) w)d\Sigma.
\end{eqnarray}
\\
{\bf Proof of (\ref{5.14}).} 
By (\ref{(5.9)}), we have the following equality:
\begin{eqnarray}\label{(5.15)}
(L_2 w, L_1 w)_{L^2 (Q)}            \nonumber\\
= - \int_Q \left( \sumij a_{ij} \ppp_i\ppp_jw - s^2\lambda^2\varphi^2 
a(x,t, \nabla\psi, \nabla\psi) w - s(\ppp_t\alpha) w\right)\nonumber\\
\times \left( \ppp_tw + 2s\lambda^2
\varphi a(x,t, \nabla\psi, \nabla\psi) w\right) dxdt\nonumber \\
-  \int_Q(2s^3\lambda^3\varphi^3  a(x,t, \nabla\psi, \nabla\psi)w
+ 2s^2 \la\va (\ppp_t\alpha) w) a(x,t, \nabla\psi, \nabla w) dxdt \nonumber\\
 -  \int_Q \left( \sumij a_{ij}\ppp_i\ppp_j w\right)
2s\lambda\varphi a(x,t, \nabla\psi, \nabla w) dxdt\nonumber\\
=:  A_1 + A_2 + A_3.
\end{eqnarray}

Now we calculate $A_1, A_2, A_3$.
\\
{\bf Calculations of $A_1$}
\\
By integrating by parts the first term on the right-hand-side, we obtain
\begin{eqnarray}\label{(5.16)}
A_1  =  \int_Q \left( -\sumij a_{ij} \ppp_i\ppp_jw -
s^2\lambda^2\varphi^2 a(x,t, \nabla\psi, \nabla\psi) w
- s(\ppp_t\alpha) w\right)\nonumber\\
\times  (\ppp_tw + 2s\lambda^2\varphi a(x,t,\nabla\psi,
\nabla\psi) w) dxdt \nonumber\\
= \int_Q \biggl\{ \partial_t w \sumij (\ppp_ja_{ij})\ppp_iw 
+ \sumij a_{ij} (\ppp_i w)(\ppp_j\ppp_t w) \nonumber \\
-  \frac{s^2\lambda^2\varphi^2}{2}a(x,t, \nabla\psi, \nabla\psi) \ppp_t (w^2)
- \frac{s\partial_t\alpha}{2}  \partial_t(w^2)
 - 2s^3\la^4 \varphi^3 a(x,t, \nabla\psi, \nabla\psi)^2 w^2  \nonumber \\
- 2s^2\lambda^2 \va (\ppp_t\alpha) a(x,t, \nabla\psi, \nabla\psi) w^2
 + 2s\lambda^2\varphi a(x,t, \nabla\psi, \nabla\psi) w 
\sumij (\ppp_ja_{ij}) \ppp_iw  \nonumber\\
+  2s\lambda^2 \varphi a(x,t, \nabla\psi, \nabla\psi) 
a(x,t, \nabla w, \nabla w)
 + 2s\lambda^2 w \sumij a_{ij} (\ppp_jw)\ppp_i
(\varphi a(x,t, \nabla\psi, \nabla\psi)) \biggr\} dx dt  \nonumber\\
- \int_\Sigma a(x,t,\nu,\nabla w)(\partial_t w + 2s\lambda^2
\varphi a(x,t, \nabla\psi, \nabla\psi) w)d\Sigma\nonumber\\\
 = \int_Q \biggl\{ \partial_t w \sumij (\partial_ja_{ij})\ppp_iw \nonumber\\
+ \frac 12\partial_t( \sumij a_{ij} (\ppp_iw)\partial_jw
- \sumij \frac{1}{2}(\partial_t a_{ij})(\ppp_iw)\ppp_jw
- \frac{s^2\lambda^2\varphi^2}{2}
a(x,t, \nabla\psi, \nabla\psi) \partial_t (\vert w\vert^2)\nonumber\\
- \frac{s\partial_t\alpha}{2}  \partial_t \vert w\vert^2
 - 2s^3 \la^4 \varphi^3 a(x,t, \nabla\psi, \nabla\psi)^2 w^2 \nonumber \\
-  2s^2\lambda^2 \va (\partial_t\alpha) a(x,t, \nabla\psi, \nabla\psi) w^2
 + 2s\lambda^2\varphi a(x,t, \nabla\psi, \nabla\psi) w 
\sumij (\ppp_ja_{ij})\ppp_jw                    \nonumber \\
+ 2s\lambda^2 \varphi a(x,t, \nabla\psi, \nabla\psi)a(x,t,\nabla w, \nabla w)
+ 2s\lambda^2 w \sumij a_{ij} (\ppp_jw)
\ppp_i(\varphi a(x,t, \nabla\psi, \nabla\psi))\biggr\} dx dt\nonumber\\    
-  \int_\Sigma a(x,t,\nu,\nabla w)(\partial_t w + 2s\lambda^2\varphi 
a(x,t, \nabla\psi, \nabla\psi) w)d\Sigma.                  
\end{eqnarray}
Integrating this equality by parts with respect to $t$, we obtain
\begin{eqnarray}\label{(5.17)}
 A_1 = 
\int_Q \biggl\{ \ppp_tw\sumij \ppp_ja_{ij} 
- \frac{s^2\lambda^2\varphi^2}{2}a(x,t, \nabla\psi, \nabla\psi) 
\partial_t (\vert w\vert^2)\nonumber\\
  + \frac{s\partial^2_t\alpha}{2}  w^2
 - 2s^3\la^4 \varphi^3 a(x,t, \nabla\psi, \nabla\psi)^2 w^2\nonumber \\
-  2s^2\lambda^2 \va (\partial_t\alpha) a(x,t, \nabla\psi, \nabla\psi) w^2
 + 2s\lambda^2\varphi a(x,t, \nabla\psi, \nabla\psi) w 
\sumij (\partial_ja_{ij}) \ppp_iw                      \nonumber      \\
+  2s\lambda^2 \varphi a(x,t, \nabla\psi, \nabla\psi) 
a(x,t, \nabla w, \nabla w)
 + 2s\lambda^2 w \sumij a_{ij} (\ppp_jw)\ppp_i
(\varphi a(x,t, \nabla\psi, \nabla\psi))\biggr\} dx dt 
\nonumber\\
-  \int_\Sigma a(x,t,\nu,\nabla w)(\partial_t w
+ 2s\lambda^2\varphi a(x,t, \nabla\psi, \nabla\psi) w)d\Sigma.                  \end{eqnarray}
\\
{\bf Calculation of $A_2$}
\\
\begin{eqnarray}\label{(5.18)}
A_2 = - \int_Q (2s^3\lambda^3\va^3 w a(x,t, \nabla\psi, \nabla\psi)
a(x,t, \nabla\psi, \nabla w)   + 2s^2 \lambda\va (\partial_t\alpha) w 
a(x,t, \nabla\psi, \nabla w)) dxdt      
\nonumber\\
=   -\int_Q (s^3\lambda^3\varphi^3 a(x,t,\nabla\psi,\nabla\psi)
a(x,t, \nabla\psi, \nabla (w^2) ) 
+ s^2\la\va (\partial_t \alpha)a(x,t, \nabla\psi, \nabla (w^2))) dxdt 
\nonumber\\
= \int_Q \biggl\{ 3s^3\lambda^4 \varphi^3 a(x,t, \nabla\psi, \nabla\psi)^2 w^2
 + s^3\la^3\va^3 w^2 \sumij \ppp_i
(a_{ij} (\ppp_j\psi) a(x,t, \nabla\psi, \nabla\psi))
\nonumber\\
 + \sumij \partial_j \left(\frac{s^2 \lambda^2\va (\ppp_t\alpha)}{2} 
a_{ij}\ppp_i\psi \right) w^2 \biggr\} 
dxdt                
\nonumber\\
- \int_\Sigma(2s^3\lambda^3 \varphi^3 a(x,t, \nabla\psi, \nabla\psi)
a(x,t, \nabla\psi, \nu) 
  + 2s^2 \lambda\va (\partial_t\alpha) a(x,t, \nabla\psi, \nu))w^2
d\Sigma.                                     
\end{eqnarray}
\\
{\bf Calculation of $A_3$}
\\
\begin{eqnarray}
 A_3  = \int_Q - \left(\sumij a_{ij}\ppp_i\ppp_j w\right)
\left( 2s\lambda\varphi \sum_{k,\ell=1}^d a_{k\ell}(\partial_k\psi)
(\partial_\ell w)\right) dxdt\nonumber \\
 = \int_Q \biggl\{
\sumij (\partial_ja_{ij})(\ppp_iw) 2s\lambda\varphi \sum_{k,\ell=1}^d 
a_{k\ell} (\partial_k\psi)(\partial_\ell w)
 + 2s\lambda^2 \varphi a(x,t, \nabla\psi, \nabla w)^2 \nonumber \\
 +  2s\lambda\varphi (\sumij a_{ij} (\ppp_iw))
\left( \sum_{k,\ell=1}^d \partial_j(a_{k\ell}\ppp_k\psi )\ppp_\ell w\right) 
 + 2s\lambda\varphi \sumij a_{ij} \partial_{x_i} w
\sum_{k,\ell=1}^d a_{k\ell}(\partial_k\psi)(\ppp_j\ppp_\ell w) \biggr\} dxdt\nonumber\\
+  \int_\Sigma 2s\lambda\varphi 
a(x,t,\nu,\nabla w)a(x,t,\nabla\psi,\nabla w) d\Sigma  \nonumber  \\
 =  \int_Q \biggl\{ 
 \sumij (\partial_j a_{ij})(\ppp_iw) 2s\lambda\varphi 
\sum_{k,\ell=1}^d a_{k\ell} (\ppp_k\psi)\partial_\ell w 
+ 2s\lambda^2 \varphi a(x,t, \nabla\psi, \nabla w)^2  \nonumber    \\
 +  2s\lambda\varphi \sumij a_{ij} (\partial_i w)
\left(  \sum_{k, \ell =1}^d \partial_j (a_{k\ell} (\partial_k\psi)
\partial_\ell w \right)
- s\lambda\varphi
\sum_{k, \ell = 1}^d a_{k\ell} (\partial_k\psi) 
\left( \sumij (\partial_\ell a_{ij})(\partial_i w)(\partial_j w)\right) \nonumber \\
 +  s\lambda \varphi \sum_{k,\ell = 1}^d a_{k\ell}(\partial_k\psi)
\partial_{x_\ell} ( \sumij a_{ij}(\ppp_iw)\ppp_jw \biggr\} dxdt\nonumber\\
+  \int_\Sigma  2s\lambda\varphi 
a(x,t,\nu,\nabla w)a(x,t,\nabla\psi,\nabla w) d\Sigma.   \nonumber      
\end{eqnarray}
Integrating by parts once again, we obtain
\begin{eqnarray}\label{(5.19)}
A_3 = \int_Q \biggl\{ \sumij (\partial_j a_{ij})(\ppp_iw) 
2s\lambda\varphi \sum_{k,\ell=1}^d a_{k\ell} (\partial_k\psi)
\ppp_{\ell}w + 2s\lambda^2 \varphi
 a(x,t, \nabla\psi, \nabla w)^2\nonumber\\
+ 2s\lambda \varphi \sumij a_{ij} (\ppp_iw) 
\left( \sum_{k,\ell=1}^d \partial_j (a_{k\ell}\partial_k\psi)
\partial_\ell w\right) 
- s\lambda\varphi \sum_{k,\ell=1}^d a_{k\ell} (\partial_k\psi)
\left( \sumij (\partial_\ell a_{ij})(\partial_i w)\partial_jw\right) 
\nonumber  \\
-  s\lambda^2 \varphi a(x,t, \nabla\psi, \nabla\psi)a(x,t, \nabla w, \nabla w)
 - s\lambda\varphi \sum_{k,\ell=1}^d a_{k\ell}(\partial_k\psi)
 \left( \sumij (\partial_\ell a_{ij})(\partial_i w) \partial_jw \right) 
\nonumber\\
 -  a(x,t, \nabla w, \nabla w) s\lambda\varphi \sum_{k,\ell=1}^d
\partial_\ell (a_{k\ell}\ppp_k\psi)\biggr\} dxdt   
\nonumber\\
+  \int_\Sigma(2s\lambda\varphi a(x,t,\nu,\nabla w)a(x,t,\nabla\psi,\nabla w) 
- s\lambda\varphi a(x,t, \nabla w, \nabla w)
a(x,t, \nu, \nabla\psi)) d\Sigma.                    
\end{eqnarray}

Taking the sum of (\ref{(5.17)}) - (\ref{(5.19)}), we obtain
\begin{eqnarray}\label{(5.20)}
(L_2 w, L_1 w)_{L^2 (Q)}
= \int_Q \biggl\{ \partial_t w \sumij (\ppp_ja_{ij})\ppp_jw 
- \sum_{j=1}^d \frac12(\partial_t a_{ij})(\ppp_iw)\ppp_jw  \nonumber  \\
+ \partial_t\left( \frac{s^2\lambda^2\varphi^2}{2}
a(x,t, \nabla\psi, \nabla\psi)\right) w^2
  +  \frac{s\partial^2_t\alpha}{2}  w^2
- 2s^3\la^4 \varphi^3 a(x,t, \nabla\psi, \nabla\psi)^2 w^2  \nonumber  \\
-  2s^2\lambda^2\va (\ppp_t\alpha) a(x,t, \nabla\psi, \nabla\psi) w^2
 + 2s\lambda^2 \varphi a(x,t, \nabla\psi, \nabla\psi) w 
\sumij (\partial_ja_{ij})(\partial_i w)   \nonumber      \\
+  2s\lambda^2 \varphi a(x,t, \nabla\psi, \nabla\psi)a(x,t,\nabla w, \nabla w)
+ 2s\lambda^2 w \sumij a_{ij} (\partial_j w)
\partial_i( \varphi a(x,t, \nabla\psi, \nabla\psi)) \nonumber  \\
+  \sumij (\partial_j a_{ij})(\ppp_iw) 2s\lambda\varphi 
\sum_{k,\ell=1}^d a_{k\ell}(\partial_k \psi) \partial_{\ell} w
 + 2s\lambda^2 \varphi a(x,t, \nabla\psi, \nabla (w^2))  \nonumber      \\
+ 2s\lambda \varphi \sumij  a_{ij} (\partial_i w)
\left( \sum_{k,\ell=1}^d \partial_j (a_{k\ell} \ppp_k\psi)\partial_\ell w
\right)
- s\lambda\varphi \sum_{k,\ell=1}^d
a_{k\ell} (\ppp_k\psi) \left(\sumij (\partial_\ell a_{ij})(\partial_i w)
\partial_j w \right)                                     \nonumber\\
 -  s\lambda^2 \varphi a(x,t, \nabla\psi, \nabla\psi) 
a(x,t, \nabla w, \nabla w)
- s\lambda\varphi \sum_{k,\ell=1}^d a_{k\ell} (\partial_k\psi)
\left( \sumij (\partial_\ell a_{ij})(\partial_i w) 
\partial_j w \right)        \nonumber           \\
- a(x,t, \nabla w, \nabla w) s\lambda\varphi \sum_{k,\ell=1}^d
\partial_\ell(a_{k\ell}\ppp_k\psi)) \biggr\} dxdt        \nonumber\\
+  \int_\Sigma (2s\lambda\varphi a(x,t,\nu,\nabla w)
a(x,t,\nabla\psi,\nabla w) 
- s\lambda\varphi a(x,t, \nabla w, \nabla w)a(x,t, \nu, \nabla\psi)) d\Sigma
                                                             \nonumber  \\
- \int_\Sigma(2s^3\lambda^3 \varphi^3 a(x,t, \nabla\psi, \nabla\psi)
a(x,t, \nabla\psi, \nu) 
+ 2s^2 \lambda\va (\partial_t\alpha) a(x,t, \nabla\psi, \nu))w^2
d\Sigma                                                                    
\nonumber\\
- \int_\Sigma a(x,t,\nu,\nabla w)(\partial_t w
+ 2s\lambda^2\varphi a(x,t, \nabla\psi, \nabla\psi) w)d\Sigma.    
\end{eqnarray}
Finally we observe that
\begin{align*}
 \partial_t w \left(\sumij (\partial_j a_{ij})\partial_i w\right)
+ 2s\lambda^2 \varphi a(x,t, \nabla\psi, \nabla\psi) w 
\sumij (\partial_j a_{ij}) \partial_i w            \\
+  2s\lambda^2 \varphi a(x,t, \nabla\psi, \nabla\psi) w 
\sumij (\partial_j a_{ij}) \partial_i w
= L_1(x,t,D) w\, \sumij (\partial_j a_{ij}) \partial_iw.
\end{align*}
Thus the proof of (\ref{5.14}) is complete.  
$\blacksquare$
\\

Similarly to (\ref{5.14}), we can readily verify 
\begin{eqnarray}\label{(5.21)}
 (\widetilde L_2\widetilde  w,\, \widetilde  L_1\widetilde  w)_{L^2 (Q)}
= \int_Q \widetilde L_1(x,t,D,s)\widetilde w 
\sumij (\partial_j a_{ij}) \partial_i\widetilde  wdxdt \nonumber\\
+ \int_Q \biggl\{ -\sumij \frac 12 (\partial_t a_{ij})(\partial_i\widetilde w)
\partial_j\widetilde w                                
+ \partial_t\left( \frac{s^2\lambda^2\varphi^2}{2}
a(x,t, \nabla\psi, \nabla\psi)\right) \widetilde w^2
  + \hhalf s(\ppp_t^2\alpha) \widetilde  w^2\nonumber\\
+ s^3\la^4 \varphi^3 a(x,t, \nabla\psi, \nabla\psi)^2\widetilde w^2
- 2s^2\lambda^2\va (\ppp_t \alpha) a(x,t, \nabla\psi, \nabla\psi) 
\widetilde w^2\nonumber\\  
+ s\lambda^2 \varphi a(x,t, \nabla\psi, \nabla\psi) 
a(x,t, \nabla\widetilde w, \nabla\widetilde  w)   \nonumber \\
 +  2s\lambda^2\widetilde w \left(\sumij a_{ij} (\partial_j\widetilde w)
\partial_i (\varphi a(x,t, \nabla\psi, \nabla\psi))\right)
+ 2s\lambda^2 \varphi a(x,t, \nabla\psi, \nabla \widetilde w)^2\nonumber\\
+  2s\lambda \varphi \sumij a_{ij} (\partial_i\widetilde w)
\left(\sum_{k,\ell=1}^d \partial_j (a_{k\ell} (\ppp_k\psi)\partial_\ell
\widetilde w \right) 
- s\lambda\varphi \sum_{k,\ell=1}^d
a_{k\ell} (\partial_k\psi) \left( \sumij (\partial_\ell a_{ij})
(\partial_i\widetilde  w) \partial_j\widetilde w \right)        \nonumber\\
 - s\lambda\varphi \sum_{k,\ell=1}^d a_{k\ell} (\ppp_k\psi)
\left( \sumij (\partial_\ell a_{ij})(\partial_i\widetilde w)
\partial_j\widetilde w \right)        
- s\la\va a(x,t, \nabla\widetilde w, \nabla\widetilde w) 
\sum_{k,\ell=1}^d \partial_\ell (a_{k\ell}\ppp_k\psi) \biggr\} dxdt \nonumber\\
 - \int_\Sigma ( 2s\lambda\varphi 
a(x,t,\nu,\nabla \widetilde w)a(x,t,\nabla\psi,\nabla \widetilde w) 
- s\va\lambda\varphi  a(x,t, \nabla\widetilde  w, \nabla \widetilde w)
a(x,t, \nu,\nabla \psi)) d\Sigma            \nonumber                                \\ 
+  \int_\Sigma(2s^3\lambda^3 \varphi^3 a(x,t, \nabla\psi, \nabla\psi)
a(x,t, \nabla\psi, \nu) + 2s^2 \lambda\va(\partial_t\alpha)
a(x,t, \nabla\psi, \nu))\widetilde w^2 d\Sigma
\nonumber\\
- \int_\Sigma a(x,t,\nu,\nabla \widetilde w)
(\partial_t \widetilde w + 2s\lambda^2 \varphi a(x,t,\nabla\psi, \nabla\psi)
\widetilde w)d\Sigma.
\end{eqnarray}
\\
{\bf Third Step: completion of the proof of Lemma \ref{6}}.
\\
Taking the $L^2-$ norms of both sides of the equations (\ref{(5.10)}) and 
(\ref{(5.12)}), 
we have
$$
\Vert L_1w\Vert^2_{L^2(Q)} + 2(L_1w,\, L_2w)_{L^2(Q)}
+ \Vert L_2w\Vert^2_{L^2(Q)} = \Vert H\Vert^2_{L^2(Q)}
$$
and 
$$
\Vert \widetilde L_1\widetilde w\Vert^2_{L^2(Q)}
+ 2(\widetilde L_1\widetilde w,\, \widetilde L_2\widetilde w)_{L^2(Q)}
+ \Vert \widetilde L_2\widetilde w\Vert^2_{L^2(Q)}
= \Vert \widetilde{H}\Vert^2_{L^2(Q)}.
$$
We take the parameter $\lambda$ sufficiently large, so that 
\begin{eqnarray}\label{(5.22)}
 (L_1 w, \, L_2 w)_{L^2 (Q)}
\ge \int_Q L_1(x,t,D,s)w \sumij (\partial_j a_{ij}) \partial_i w dxdt   
                                                   \nonumber   \\
+  \frac 14 \int_Q( s\lambda^2 \varphi
a(x,t, \nabla\psi, \nabla\psi) a(x,t, \nabla w, \nabla w)
+ s^3 \la^4\varphi^3 a(x,t, \nabla\psi, \nabla\psi)^2 w^2) dxdt  
                                           \nonumber      \\
+   \int_\Sigma ( 2s\lambda\varphi a(x,t,\nu,\nabla w)
a(x,t,\nabla\psi,\nabla w) 
- s\lambda\varphi a(x,t, \nabla w, \nabla w) a(x,t, \nu,\nabla \psi)) d\Sigma
                                                              \nonumber  \\ 
-  \int_\Sigma (2s^3\lambda^3 \varphi^3 a(x,t, \nabla\psi, \nabla\psi)
a(x,t, \nabla\psi, \nu)   
+ 2s^2 \lambda\va (\partial_t\alpha) a(x,t, \nabla\psi, \nu))w^2 d\Sigma  
                                                       \nonumber  \\
- \int_\Sigma a(x,t,\nu,\nabla w)(\partial_t w 
+ 2s\lambda^2\varphi a(x,t, \nabla\psi, \nabla\psi) w)d\Sigma.
                                      \end{eqnarray}
Using (\ref{5.14}), we have 
\begin{eqnarray}\label{(5.23)}
 \sum_{k=1}^2\Vert L_kw\Vert^2_{L^2(Q)}
+ 2\int_Q L_1(x,t,D,s)w \sumij (\partial_j a_{ij})\partial_iw dxdt  
                                                      \nonumber   \\
+ \frac 12 \int_Q( s\lambda^2 \varphi
a(x,t, \nabla\psi, \nabla\psi) a(x,t, \nabla w, \nabla w)
+ s\la^4\varphi^3 a(x,t, \nabla\psi, \nabla\psi)^2 w^2)dxdt   \nonumber\\
+  \int_\Sigma (2s\lambda\varphi a(x,t,\nu,\nabla w)a(x,t,\nabla\psi,\nabla w)
- s\lambda\varphi a(x,t, \nabla w, \nabla w)
a(x,t, \nu,\nabla \psi)) d\Sigma     \nonumber\\ 
-  \int_\Sigma(2s^3\lambda^3 \varphi^3 a(x,t, \nabla\psi, \nabla\psi)
a(x,t, \nabla\psi, \nu) + 2s^2 \lambda \va (\partial_t\alpha)
a(x,t, \nabla\psi, \nu))w^2 d\Sigma                                         
\nonumber\\
- \int_\Sigma a(x,t,\nu,\nabla w)(\partial_t w
+ 2s\lambda^2\varphi a(x,t, \nabla\psi, \nabla\psi) w)d\Sigma
\le \Vert H \Vert^2_{L^2(Q)}.                       
\end{eqnarray}
This inequality implies 
\begin{eqnarray}\label{(5.24)}
 \frac 14 \sum_{k=1}^2\Vert L_kw\Vert^2_{L^2(Q)}\nonumber\\
+ \frac 14 \int_Q ( s\lambda^2 \varphi a(x,t, \nabla\psi, \nabla\psi) 
a(x,t, \nabla w, \nabla w)
+ s^3 \la^4\varphi^3 a(x,t, \nabla\psi, \nabla\psi)^2 w^2)dxdt    \nonumber  \\
+  \int_\Sigma(2s\lambda\varphi a(x,t,\nu,\nabla w)a(x,t,\nabla\psi,\nabla w)
 - s\lambda\varphi a(x,t, \nabla w, \nabla w)a(x,t, \nu,\nabla \psi)) d\Sigma
                                                  \nonumber\\ 
-  \int_\Sigma(2s^3\lambda^3\varphi^3 a(x,t, \nabla\psi, \nabla\psi)
a(x,t, \nabla\psi,\nu)
+ 2s^2 \lambda\va (\partial_t\alpha) a(x,t, \nabla\psi, \nu))w^2 d\Sigma     
\nonumber\\
- \int_\Sigma a(x,t,\nu,\nabla w)(\partial_t w
+ 2s\lambda^2 \varphi a(x,t, \nabla\psi, \nabla\psi) w)d\Sigma
\le \Vert H \Vert^2_{L^2(Q)}.
\end{eqnarray}
Similarly, using (\ref{(5.21)}), we obtain
\begin{eqnarray}\label{(5.25)}
 \frac 14 \sum_{k=1}^2\Vert \widetilde L_k\widetilde w\Vert^2_{L^2(Q)} 
                                                            \nonumber \\
+  \frac 14 \int_Q( s\lambda^2 \varphi a(x,t, \nabla\psi, \nabla\psi) 
a(x,t, \nabla \widetilde w, \nabla\widetilde  w)
+ s^3 \la^4\varphi^3 a(x,t, \nabla\psi, \nabla\psi)^2\widetilde w^2)dxdt
                                                              \nonumber \\
-  \int_\Sigma (2s\lambda\varphi a(x,t,\nu,\nabla\widetilde  w)
a(x,t,\nabla\psi,\nabla\widetilde  w) 
- s\lambda\varphi a(x,t, \nabla\widetilde w, \nabla\widetilde  w)
a(x,t, \nu,\nabla \psi)) d\Sigma                            \nonumber\\ 
+  \int_\Sigma(2s^3\lambda^3\varphi^3 a(x,t, \nabla\psi, \nabla\psi)
a(x,t, \nabla\psi, \nu) 
+ 2s^2 \lambda\va(\partial_t\alpha) a(x,t, \nabla\psi, \nu))\widetilde w^2
d\Sigma
\nonumber\\
-\int_\Sigma a(x,t,\nu,\nabla\widetilde  w)(\partial_t\widetilde  w
+ 2s\lambda^2 \varphi a(x,t, \nabla\psi, \nabla\psi) \widetilde w)d\Sigma
\le \Vert \www{H}\Vert^2_{L^2(Q)}.
\end{eqnarray}
We set 
\begin{eqnarray}\label{(5.26)}
I := \int_{(\ppp\OOO\setminus \Gamma)\times (0,T)} 
(2s\lambda\varphi a(x,t,\nu,\nabla w)
a(x,t,\nabla\psi,\nabla w) 
- s\lambda\varphi a(x,t, \nabla w, \nabla w)
a(x,t, \nu,\nabla \psi)) d\Sigma        \nonumber\\ 
- \int_{(\ppp\OOO\setminus \Gamma)\times (0,T)}
(2s^3\lambda^3\varphi^3 a(x,t, \nabla\psi, \nabla\psi)
a(x,t, \nabla\psi, \nu) + 2s^2 \lambda\va (\partial_t\alpha)
a(x,t, \nabla\psi, \nu))w^2 d\Sigma                      \nonumber\\
-  \int_{(\ppp\OOO\setminus \Gamma)\times (0,T)} 
a(x,t,\nu,\nabla w)(\partial_t w
+ 2s\lambda^2\varphi a(x,t, \nabla\psi, \nabla\psi) w)d\Sigma  \nonumber \\
- \int_{\Sigma} 
(2s\lambda\varphi a(x,t,\nu,\nabla\widetilde w)
a(x,t,\nabla\psi,\nabla\widetilde w) 
- s\lambda\varphi a(x,t, \nabla \widetilde w, \nabla\widetilde w)
a(x,t, \nu,\nabla \psi)) d\Sigma                \nonumber\\ 
+  \int_{\Sigma}
(2s^3\lambda^3\varphi^3 a(x,t, \nabla\psi, \nabla\psi)
a(x,t, \nabla\psi, \nu) + 2s^2 \lambda\va (\partial_t\alpha)
a(x,t, \nabla\psi, \nu))\widetilde w^2 d\Sigma 
\nonumber\\
- \int_{\Sigma} a(x,t, \nu,\nabla\widetilde w)
(\partial_t\widetilde w + 2s\lambda^2\varphi a(x,t, \nabla\psi, \nabla\psi) 
\widetilde w)d\Sigma.
\end{eqnarray}
We observe that
$$
\nabla w= (\nabla u)e^{s\alpha} + s\lambda \varphi (\nabla \psi) ue^{s\alpha}, 
\quad 
\nabla \widetilde w= (\nabla u)e^{s\alpha} - s\lambda \varphi (\nabla \psi)
u e^{s\alpha}\quad \mbox{on}\quad (\ppp\OOO \setminus \Gamma) \times (0,T).
$$
Therefore, we rewrite (\ref{(5.26)}) as 
\begin{eqnarray}\label{(5.27)}
 I = \int_{(\ppp\OOO\setminus \Gamma)\times (0,T)} 
\{ s^2\lambda^2\varphi^2 
(2a(x,t,\nu,\nabla u)a(x,t,\nabla\psi, \nabla \psi)
+ 2a(x,t,\nu,\nabla \psi)a(x,t,\nabla\psi,\nabla u)\nonumber \\
- 2s a(x,t,\nabla\psi,\nabla u) a(x,t, \nu,\nabla \psi)u\} \weight d\Sigma
                                                       \nonumber   \\ 
- \int_{(\ppp\OOO\setminus \Gamma)\times (0,T)}
a(x,\,t,\,\nu,\,(\nabla u)e^{s\alpha}+s\lambda\varphi(\nabla \psi)ue^{s\alpha})
\,\,
(\partial_t w + 2s\lambda^2\varphi a(x,t, \nabla\psi, \nabla\psi) w)d\Sigma
                                                                  \nonumber \\
- \int_{(\ppp\OOO\setminus \Gamma)\times (0,T)} 
a(x,\,t,\,\nu,(\nabla u)e^{s\alpha}-s\lambda\varphi(\nabla \psi)ue^{s\alpha})
\,\,(\partial_t w + 2s\lambda^2\varphi a(x,t, \nabla\psi, \nabla\psi) w)d\Sigma
                                         \nonumber       \\
=  \int_{(\ppp\OOO\setminus \Gamma)\times (0,T)} 
2s^2\lambda^2\varphi^2 a(x,t,\nu,\nabla u)a(x,t,\nabla\psi, \nabla \psi)
ue^{2s\alpha} d\Sigma
\nonumber\\
- 2\int_{(\ppp\OOO\setminus \Gamma)\times (0,T)} 
a(x,t,\nu,(\nabla u)e^{s\alpha})(\partial_t w
+ 2s\lambda^2\varphi a(x,t, \nabla\psi, \nabla\psi) w)d\Sigma.
\end{eqnarray}
Since the Robin boundary condition implies
$$
\partial_{\nu_A}u-p(x,t)u = g \quad \mbox{on $(\ppp\OOO\setminus \Gamma)
\times (0,T)$}
$$
then, by (\ref{(5.27)}), we obtain
\begin{eqnarray}
 I = \int_{(\ppp\OOO\setminus \Gamma)\times (0,T)} 
s^2\lambda^2\varphi^2 
2a(x,t,\nabla \psi,\nabla \psi)(p(x,t)u+g)ue^{2s\alpha} d\Sigma \\ 
- 2\int_{(\ppp\OOO\setminus \Gamma)\times (0,T)} (p(x,t)u+g)e^{s\alpha}
(\partial_t w + 2s\lambda^2 \varphi a(x,t, \nabla\psi,\nabla\psi)w)d\Sigma
                                            \nonumber \\
= \int_{(\ppp\OOO\setminus \Gamma)\times (0,T)} 
 2s^2\lambda^2\varphi^2 a(x,t,\nabla \psi,\nabla \psi)(p(x,t)u+g)ue^{2s\alpha} 
d\Sigma  \nonumber\\
+ 2\int_{(\ppp\OOO\setminus \Gamma)\times (0,T)}
\left( (\partial_tp)(x,t)\frac{w^2}{2}+ \partial_t (ge^{s\alpha})w
- 2s\lambda^2 \varphi a(x,t, \nabla\psi, \nabla\psi) 
(p(x,t)w+ge^{s\alpha})w \right)d\Sigma.\nonumber
\end{eqnarray}
Hence,  for any positive $\epsilon$ we have
\begin{equation}\label{(5.28)}
\vert I\vert\le \int_{(\ppp\OOO\setminus \Gamma)\times (0,T)}
\left(\epsilon s^\frac 52\lambda^2\varphi^\frac 52w^2 + C(\epsilon)\left(
\frac{\vert \partial_t g\vert^2e^{2s\alpha}}
{s^2\lambda^2\varphi^2} + \frac{1}{\root\of{s\varphi}} g^2e^{2s\alpha}\right)
\right)d\Sigma.
\end{equation}
By (\ref{(5.24)}), (\ref{(5.25)}) and (\ref{(5.28)}), we have 
\begin{eqnarray}\label{(5.29)}
 \frac 14 \sum_{k=1}^2\Vert L_kw\Vert^2_{L^2(Q)} 
+ \frac 14 \sum_{k=1}^2\Vert \widetilde L_k\widetilde w\Vert^2_{L^2(Q)}
                                              \nonumber\\
+ \frac 14 \int_Q( s\lambda^2 \varphi
a(x,t, \nabla\psi, \nabla\psi) a(x,t, \nabla w, \nabla w)
+ s^3 \la^4\varphi^3 a(x,t, \nabla\psi, \nabla\psi)^2 w^2)dxdt \nonumber\\
+  \int_{\Gamma \times (0,T)}(2s\lambda\varphi a(x,t,\nu,\nabla w)
a(x,t,\nabla\psi,\nabla w) 
- s\lambda\varphi a(x,t, \nabla w, \nabla w)a(x,t, \nu,\nabla \psi))d\Sigma
                                                         \nonumber      \\ 
-  \int_{\Gamma \times (0,T)}(2s^3\lambda^3\varphi^3 
a(x,t, \nabla\psi, \nabla\psi)a(x,t, \nabla\psi, \nu)   
+ 2s^2 \lambda\va (\partial_t\alpha) a(x,t, \nabla\psi, \nu))w^2 d\Sigma 
                                                              \nonumber\\
-  \int_{\Gamma \times (0,T)} a(x,t,\nu,\nabla w)(\partial_t w
+ 2s\lambda^2\varphi a(x,t, \nabla\psi, \nabla\psi) w)d\Sigma   \nonumber  \\ 
 -   \int_{\Gamma \times (0,T)} 
(2s\lambda\varphi a(x,t,\nu,\nabla\widetilde w)
a(x,t,\nabla\psi,\nabla\widetilde  w) 
 -s\lambda\varphi a(x,t, \nabla\widetilde w, \nabla\widetilde  w)
a(x,t, \nu,\nabla \psi)) d\Sigma                           \nonumber\\ 
+  \int_{\Gamma \times (0,T)}(2s^3\lambda^3\varphi^3 
a(x,t,\nabla\psi, \nabla\psi)a(x,t, \nabla\psi, \nu)   
+ 2s^2 \lambda\va(\partial_t\alpha)a(x,t, \nabla\psi, \nu))\widetilde w^2
d\Sigma                      \nonumber   \\
-  \int_{\Gamma \times (0,T)} a(x,t,\nu,\nabla\widetilde w)
(\partial_t\widetilde w + 2s\lambda^2 \varphi a(x,t, \nabla\psi, \nabla\psi) 
\widetilde w)d\Sigma
\le \Vert\widetilde{H}\Vert^2_{L^2(Q)}
+ \Vert H \Vert^2_{L^2(Q)}\nonumber\\
+ \int_{(\ppp\OOO\setminus \Gamma)\times (0,T)}
\left(\epsilon s^\frac 52\lambda^2\varphi^\frac 52w^2 + C(\epsilon)\left(
\frac{\vert \partial_t g\vert^2e^{2s\alpha}}
{s^2\lambda^2\varphi^2} + \frac{1}{\root\of{s\varphi}} g^2e^{2s\alpha}\right)
\right)d\Sigma.                        
\end{eqnarray}
We note  
\begin{equation}\label{(5.30)}
\Vert\widetilde{H} \Vert^2_{L^2(Q)}
\le C\Vert H \Vert^2_{L^2(Q)}.                   
\end{equation}

Taking the scalar product of  functions $L_1(x,t,D,s) w$ and 
$s^\frac 32\lambda^2 \varphi^\frac 32 w$ in $L^2(Q)$ and integrating by parts 
we obtain
\begin{align*}
& (L_1(x,t,D,s) w,s^\frac 32\lambda^2 \varphi^\frac 32 w)_{L^2(Q)}\\
= & (\partial_t w + 2{{s}}\lambda \varphi
\sum_{i,j=1}^n a_{ij} \psi_{x_i} \partial_{x_j} w + 2{{s}}\lambda^2 \varphi
a(t, x, \nabla\psi, \nabla\psi) w,s^\frac 32\lambda^2 \varphi^\frac 32 w)
_{L^2(Q)}\\
=& \int_Q \left(-\frac 12 \partial_t (s^\frac 32 \lambda^2\varphi^\frac 32)
w^2 - \sum_{i,j=1}^n s\lambda\partial_{x_j}(\varphi a_{ij}\psi_{x_i})w^2
+ 2s^\frac 52 \lambda^4 \varphi^\frac 52 w^2\right)dx \\
+& \int_\Sigma s^\frac 52\lambda^2 \varphi^\frac 52 a(x,t,\nu,\nabla \psi) 
w^2d\Sigma.
\end{align*}

The above equality implies
\begin{eqnarray}\label{(5.31)}
\int_{(\ppp\OOO\setminus \Gamma)\times (0,T)} 
s^\frac 52\lambda^2 \varphi^\frac 52 w^2d\Sigma
                                                \nonumber\\
\le C\int_Q(s\lambda^2 \varphi a(x,t, \nabla\psi, \nabla\psi) 
a(x,t, \nabla w, \nabla w)
+ s^3 \varphi^3 \lambda^4 a(x,t, \nabla\psi, \nabla\psi)^2 w^2)dxdt
                                                            \nonumber\\
+ C\int_{\Gamma \times (0,T)}s^\frac 52\lambda^2 \varphi^\frac 52  w^2d\Sigma.
                                                            \nonumber 
\end{eqnarray}

By (\ref{(5.29)}) - (\ref{(5.31)}), we obtain 
\begin{eqnarray}\label{(5.32)}
\int_Q ( s\va \vert \nabla u\vert^2 + s^3 \varphi^3 \vert u\vert^2)
e^{2s\alpha} dx\,dt
+ \sum_{k=1}^2\Vert L_k(x,t,D,s) (ue^{s\alpha})\Vert^2_{L^2(Q)}  
\nonumber\\
\le C\int_Q |F|^2 e^{2s\alpha} dxdt 
+ C\int_{(\ppp\OOO\setminus \Gamma)\times (0,T)}
\left( \frac{\vert \partial_t g\vert^2e^{2s\alpha}}{s^2\lambda^2\varphi^2}
+ \frac{1}{\root\of{s\varphi}} \vert g\vert^2e^{2s\alpha} \right)d\Sigma
\nonumber\\
+ C\int_{\Gamma\times (0,T)} \left(
s\lambda\varphi \vert \nabla u\vert^2 + s^3\lambda^3\varphi^3 \vert u\vert^2 
+ \frac{\vert \partial_t u\vert^2}{s\varphi}\right) e^{2s\alpha}  d\Sigma.
\end{eqnarray}
From the definition of the operator $L_2$ and (\ref{(5.32)}), we have
\begin{eqnarray}\label{(5.33)}
 \int_Q \frac{1}{s\varphi} \vert\partial_t u\vert^2e^{2s\alpha}dxdt 
\le C\left(\Vert L_1w\Vert^2_{L^2(Q)} +
\int_Q (s\varphi |\nabla u|^2 + s^3 \varphi^3 u^2) \weight dxdt\right) 
                                            \nonumber\\
\le C\int_Q |F|^2 \weight dxdt 
+ C\int_{(\ppp\OOO\setminus \Gamma)\times (0,T)}
\left( \frac{\vert \partial_t g\vert^2\weight}{s^2\lambda^2\varphi^2}
+ \frac{1}{\root\of{s\varphi}} \vert g\vert^2e^{2s\alpha}\right)d\Sigma
                                      \nonumber\\
+ C\int_{\Gamma \times (0,T)} \left(s\lambda\varphi \vert \nabla u\vert^2 
+ s^3\lambda^3\varphi^3 \vert u\vert^2 
+ \frac{\vert \partial_t u\vert^2}{s\varphi} \right)\weight d\Sigma.
\end{eqnarray}
On the other hand 
\begin{eqnarray}\label{(5.34)}
\int_Q \frac{1}{s\varphi} \sumij \vert\ppp_i\ppp_j u\vert^2e^{2s\alpha}dxdt
\le C\biggl( \Vert L_1w\Vert^2_{L^2(Q)}
+ \int_Q (s\lambda^2\varphi |\nabla u|^2 + s^3 \lambda^3\varphi^3 u^2)
\weight dxdt                                \nonumber\\
+  \Vert ge^{s\alpha}\Vert^2_{L^2(0,T;H^\frac 12(\partial\Omega))}
+ \Vert s\varphi w\Vert^2_{L^2((\ppp\OOO\setminus \Gamma)\setminus (0,T))}
\biggr).
\end{eqnarray}
From (\ref{(5.32)}) - (\ref{(5.34)}), we complete the proof of Lemma \ref{6}. 
$\blacksquare$
\\
\vspace{0.2cm}

{\bf Acknowledgments.}
The third author was supported partly  
by Grant-in-Aid for Scientific Research (A) 20H00117 
and Grant-in-Aid for Challenging Research (Pioneering) 21K18142 of 
Japan Society for the Promotion of Science.


\begin{thebibliography}{99} %

\bibitem{ACDPS}
Y. Achdou, P. Cardaliaguet, F. Delarue, A. Porretta and F. Santambrogio, 
{\it Mean field games},
Cetraro, Italy 2019, Lecture Notes in Mathematics, C.I.M.E. 
Foundation Subseries, Volume 2281, Springer, 2019.

\bibitem{Ad}
R.A. Adams. {\it Sobolev Spaces}, Academic Press, New York, 1975.

\bibitem{CIK} D. Chae, O. Imanuvilov and S. M.  Kim, 
{\it Exact controllability for semilinear parabolic equations with Neumann 
boundary conditions,} J. Dynam. Control Systems {\bf 2} (1996) 449-483.

\bibitem{BK}
A.L. Bukhgeim and M.V. Klibanov, 
{\it Global uniqueness of a class of multidimensional inverse problems,}
Sov. Math.-Dokl. {\bf 24} (1981) 244-247.

\bibitem{CCP}
P. Cardaliaguet, M. Cirant and A. Porretta, {\it
Splitting methods and short time existence for the master equations in mean 
field games,} J. Eur. Math. Soc. (2022),
DOI 10.4171/JEMS/1227

\bibitem{FI}
A. V. Fursikov and O. Y. Imanuvilov,  {\it Controllability of Evolution 
Equations}, Lecture Notes Series vol 34, 1996, Seoul National University.

\bibitem{Hor} L. H\"ormander, {\it Linear Partial Differential Operators},
Springer-Verlag, Berlin, 1976.

\bibitem{Ima}
O.Y. Imanuvilov, {\it Controllability of parabolic equations,} Sbornik Math. 
{\bf 186} (1995) 879-900.

\bibitem{ILY}
H. Liu, O.Y. Imanuvilov and M. Yamamoto,
{\it Unique continuation for a mean field game system,}
Appl. Math. Letters 
Volume 145, November 2023, 108757 \\
https://doi.org/10.1016/j.aml.2023.108757

\bibitem{IY98}
O.Y. Imanuvilov and M. Yamamoto,
{\it Lipschitz stability in inverse parabolic problems by the Carleman
estimate,} Inverse Problems {\bf 14} (1998) 1229-1245.

\bibitem{Kli}
M.V. Klibanov,  {\it Inverse problems and Carleman estimates,} 
Inverse Problems {\bf 8} (1992) 575-596.

\bibitem{Kl23}
M. V. Klibanov, 
{\it The mean field games system: Carleman estimates, Lipschitz stability and 
uniqueness,} preprint arXiv:2303.03928

\bibitem{KlAv}
M. V. Klibanov and Y. Averboukh,
{\it Lipschitz stability estimate and uniqueness in the retrospective analysis 
for the mean field games system via two Carleman estimates,}
preprint arXiv:2302.10709


\bibitem{KLL1}
M. V. Klibanov, J. Li and H. Liu,
{\it H\"older stability and uniqueness for the mean field games system 
via Carleman estimates,} preprint, arXiv:2304.00646

\bibitem{KLL2}
M. V. Klibanov, J. Li and H. Liu,
{\it On the mean field games system with the lateral Cauchy data via 
Carleman estimates,} preprint, arXiv:2303.07556 

\bibitem{LL}
J.-M. Lasry and P.-L. Lions, {\it Mean field games, }
Japanese Journal of Mathematics, {\bf 2} (2007) 229-260.

\bibitem{LM}
J.-L. Lions and E. Magenes, {\it
Non-homogeneous Boundary Value Problems and Applications},
vol.1, Springer-Verlag, Berlin, 1972.

\bibitem{LY}
H. Liu and M. Yamamoto, {\it Stability in determination of states for the 
mean field game equations,} preprint.

\bibitem{LZ1}
H. Liu and S. Zhang, {\it On an inverse boundary problem for mean field games,}
preprint, arXiv:2212.09110 

\bibitem{LZ2}
H. Liu and S. Zhang, {\it Simultaneously recovering running cost and 
Hamiltonian 
in mean field games system,} preprint,
arXiv:2303.13096


\bibitem{Y09}
M. Yamamoto, {\it Carleman estimates for parabolic equations 
and applications,} 
Inverse Problems {\bf 25} (2009) 123013
\end{thebibliography}
\end{document}